\documentclass[11pt]{article}
\usepackage[nohead,margin=1.0in]{geometry}
\usepackage{amssymb, amsmath, amsthm,amsfonts}
\usepackage{graphicx,epsfig}
\usepackage[tight]{subfigure}
\usepackage{cite}
\usepackage{times}
\usepackage{float}
\usepackage{appendix}
\usepackage{color,xcolor}
\usepackage{epstopdf}
\usepackage{mathrsfs}
\usepackage{algorithm,algorithmic}
\usepackage{booktabs}
\usepackage{verbatim}
\usepackage{enumitem}
\usepackage[colorlinks,linktocpage,linkcolor=blue]{hyperref}
\numberwithin{equation}{section}

\newtheorem{theorem}{Theorem}[section]
\newtheorem{assumption}{Assumption}[section]
\newtheorem{lemma}{Lemma}[section]
\newtheorem{proposition}{Proposition}[section]
\newtheorem{corollary}{Corollary}[section]

\theoremstyle{remark}
\newtheorem{remark}[theorem]{Remark}

\def\d{\mathrm{d}}

\newcommand{\R}{\mathbb{R}}
\newcommand{\PP}{\mathbb{P}}
\newcommand{\tTH}{\mathcal{T}_H}
\newcommand{\FH}{\mathcal{F}_H^{\Gamma_0}}
\newcommand{\norm}[2]{\left\|#1\right\|_{#2}}
\newcommand{\seminorm}[2]{\left|#1\right|_{#2}}
\newcommand{\ip}[3]{\left(#1,#2\right)_{#3}}

\providecommand{\Tr}{\operatorname{Tr}}
\providecommand{\diam}{\operatorname{diam}}

\title{Diffusion Coefficient Identification in Source-Free Subdiffusion from a Single Terminal Observation\thanks{The work of B. Jin is supported by Hong Kong RGC General Research Fund (Project 14306824), Hong Kong RGC  ANR / RGC Joint
Research Scheme (A-CUHK402/24). The work of Z. Zhou is supported by National Natural Science Foundation of China (Project 12422117) and Hong Kong
Research Grants Council (15306525).}}
\author{Fuchang Huo\thanks{Department of Applied Mathematics, The Hong Kong Polytechnic University, Kowloon, Hong Kong,  P.R. China. (\texttt{fuchang.huo@polyu.edu.hk, zhizhou@polyu.edu.hk})}
\and Bangti Jin\thanks{Department of Mathematics, The Chinese University of Hong Kong, Shatin, New Territories, Hong Kong, P.R. China. (\texttt{b.jin@cuhk.edu.hk, bangti.jin@gmail.com})} 
\and Zhi Zhou\footnotemark[2]}
\date{}

\begin{document}

\maketitle

\begin{abstract}
In this work, we investigate the recovery of a spatially dependent diffusion coefficient in a source-free subdiffusion equation from one single terminal observation. The main challenge in the source-free setting lies in the fact that the terminal state converges to zero and thus provides no nontrivial steady-state profile. 
To address the challenge, we decompose a rescaled terminal state into a leading term characterized by an elliptic equation and a remainder: the former involves the unknown diffusion coefficient, and the latter can be absorbed into a coercive estimate associated with the leading term for large time. We establish a conditional H\"older stability estimate under asymmetric coefficient regularity: the reference coefficient is Lipschitz continuous, whereas the comparison coefficient is only in $H^1(\Omega)$.
Moreover, by suitably adapting the stability argument and nonsmsooth data analysis of the subdiffusion model, we derive error bounds for a fully discrete regularized least-squares method using finite elements in space and backward Euler convolution quadrature in time.  
Further, we investigate the joint recovery of the diffusion coefficient and the fractional order when the diffusion coefficient is known near a boundary patch, establish a boundary trace identity that determines the fractional order without requiring knowledge of the coefficient in the interior. Then we develop a two-stage reconstruction procedure in which the fractional order is recovered first, followed by the diffusion coefficient. We also derive error estimates for the order and the coefficient for the fully discrete reconstruction.
\vskip5pt

\end{abstract}

\section{Introduction}
In this work, we investigate the inverse problem of recovering a spatially varying diffusion coefficient in a subdiffusion equation from a single terminal observation. Let $\Omega\subset\mathbb{R}^d$, $d=1,2,3$, be an open bounded domain
with a $C^{1,1}$ boundary. The subdiffusion model then reads
\begin{equation}\label{eqn:fde}
    \left\{
    \begin{aligned}
        \partial_t^\alpha u-\nabla\cdot(q\nabla u) &= 0,
        && \text{in } \Omega\times(0,T),\\
        u &= 0,
        && \text{on } \partial\Omega\times(0,T),\\
        u(0) &= u_0,
        && \text{in } \Omega,
    \end{aligned}
    \right.
\end{equation}
where $T>0$ is the final time, $q$ is the unknown diffusion coefficient, and $u_0$ is the initial condition. The operator $\partial_t^\alpha$ denotes the left-sided Caputo fractional derivative of order $\alpha\in(0,1)$ in time, defined by \cite[p. 92]{kilbas2006theory}
\[
    \partial_t^\alpha u(t)
    :=\frac{1}{\Gamma(1-\alpha)}
    \int_0^t(t-s)^{-\alpha}\partial_su(s)\,\mathrm{d}s,
\]
where
$\Gamma(z)=\int_0^\infty e^{-s}s^{z-1}\,\mathrm{d}s$, for $\Re(z)>0$, is Euler's Gamma function. The diffusion coefficient $q$ belongs to the admissible set
\begin{equation*}
    \mathcal{A}:=\left\{q\in H^1(\Omega): c_0\le q(x)\le c_1 \text{ a.e. in }\Omega\right\},
\end{equation*}
where $0<c_0<c_1<\infty$ are fixed constants. Below we denote by $u(q)$ the solution of problem \eqref{eqn:fde}, to indicate its dependence on the coefficient $q$.

The model \eqref{eqn:fde} is frequently employed to describe subdiffusive transport processes, characterized by a mean squared displacement that grows sublinearly in time. Such models have attracted considerable interest in physics, biology, and finance. The list of successful applications includes subsurface flow \cite{giona1992fractional,metzler2000random}, thermal diffusion in media with fractal geometry \cite{nigmatullin1986realization}, transport in column experiments \cite{hatano1998dispersive}, and transport in highly heterogeneous aquifers \cite{adams1992field}. In the limit $\alpha\to1^-$, problem \eqref{eqn:fde} reduces to the classical diffusion model. Note that in the subdiffusive regime $0<\alpha<1$, however, solutions exhibit markedly different regularity and long-time behavior, which present substantial challenges for both mathematical analysis and numerical approximation. During the last two decades, these distinctive features have motivated intensive research on time-fractional evolution equations, including well-posedness and regularity \cite{Allen:2016,Sakamoto:2011}, long-time decay \cite{VergaraZacher:2015}, and error estimates for nonsmooth solutions \cite{JinLazarovZhou:2013,Stynes:2017,JinLiZhou:2017}. We refer interested readers to the monographs on regularity theory \cite{Kubica:book2020,Jin:book2021}, inverse problems \cite{KaltenbacherRundell:book2023}, and numerical analysis \cite{JinZhou:2023book} for further details.

In this work, we investigate the inverse problem of recovering the unknown diffusion coefficient $q^\dagger\in\mathcal{A}$ in the model \eqref{eqn:fde} from one noisy terminal observation 
\begin{equation*}
    z^\delta=z^\dag +\xi,
    \quad \text{in } \Omega, \quad \mbox{with }z^\dag=  u(q^\dagger)(\cdot,T),
\end{equation*}
where $\xi\in L^2(\Omega)$ denotes the data noise,  with a noise level $\delta:=\|z^\delta-z^\dagger\|_{L^2(\Omega)}$. The inherent ill-posedness of the inverse problem and the nonlinear dependence of the terminal state $u(q)(T)$ on the coefficient $q$ pose substantial challenges for stable reconstruction. In this work, we establish a conditional stability estimate and then analyze a regularized numerical method in two cases:
(i) the order $\alpha$ is known, and the  coefficient $q$ is unknown;
(ii) both coefficient $q$ and order $\alpha$ are unknown, under additional regularity on the initial datum $u_0$, and agreement with a prescribed coefficient near a portion $\Gamma\subset \partial\Omega$.
%These assumptions allow us to identify the order through a boundary trace identity and then recover the interior coefficient.

\subsection{Existing results}
Diffusion coefficient identification for  standard parabolic problems (i.e., $\alpha=1$) from terminal data was studied in \cite{AlessandriniVessella:1985,Isakov:1991}. Alessandrini and Vessella \cite{AlessandriniVessella:1985} investigated the inverse problem with a nonzero source in the two-dimensional case, motivated by identification of transmissivity $q$ in groundwater flow
by piezometric head measurements $u$ in  hydrogeology. Isakov \cite[Corollary 4.2]{Isakov:1991} proved the unique identifiability on an interval, with zero initial condition, zero Neumann boundary condition at one end point and a nonzero Dirichlet boundary condition at the other end point. However, in the general case, the inverse problem is still not fully understood. Indeed, in 2020, Alessandrini \cite[Section 8]{Alessandrini:2020} listed it as one of few notable open problems for inverse problems for PDEs, and distinguished two closely related formulations: source driven (nonzero source, zero initial and boundary conditions) and source free (zero source, and nonzero initial and boundary condition). There have been notable recent advances on the inverse diffusivity problem in the source free setting of the classical diffusion equation \cite{Triki:2021,LinLin:2026}.

Triki \cite{Triki:2021} developed a novel approach by exploiting the separation of exponential decay rates among eigenmodes.
Let $\lambda_j(q)$ denote the $j$th eigenvalue of the elliptic operator $-\nabla\cdot(q\nabla\cdot)$ with a zero Dirichlet boundary condition. Then the $j$th eigenmode decays as $e^{-\lambda_j(q)T}$. If the initial datum $u_0$ has a nonzero projection onto the first eigenfunction, the
first eigenmode dominates at large times, while the higher modes decay exponentially faster. Under suitable assumptions on the admissible coefficient $q$, spectral gaps and initial datum $u_0$, Triki \cite[Theorem 1.1]{Triki:2021}  established the following Lipschitz stability estimate
\begin{equation}\label{eqn:Triki}
    \|q_1-q_2\|_{L^2(\Omega)} \le c_T \|u(q_1)(T)-u(q_2)(T)\|_{H^2(\Omega)},
\end{equation}
for sufficiently large $T$, with $c_T$ growing exponentially in $T$. The analysis decomposes $-\partial_tu(q)(T)$ into $\lambda_1(q)u(q)(T)$ and a remainder involving only higher eigenmodes. Spectral perturbation estimates control the difference of these remainders for two coefficients. The faster decay provides the relative smallness needed to absorb the remainder despite the decay of the terminal state. 

Very recently, Lin and Lin \cite{LinLin:2026} investigated the standard parabolic problem in the source free setting with an exponentially growing Dirichlet input $e^{\lambda t}$ and initial datum $u_0\equiv 1$, and proved Lipschitz stability at any prescribed positive
observation time $T$. Specifically, for any $q_1,q_2\in C^{1,\gamma}(\overline{\Omega})$, $\gamma\in(0,1)$, with an a priori bound $M$, then for every $T_*>0$, there exist $\lambda$ and $C_*$ such that
\begin{align*}
  \|q_1-q_2\|_{L^2(\Omega)} \leq c\|u(q_1)(T)-u(q_2)(T)\|_{H^2(\Omega)},\quad \forall T\geq T_*.  
\end{align*}
The analysis relies crucially on a nontrivial stationary elliptic
profile obtained after suitable exponential normalization, and the terminal data problem can be viewed as a perturbation of the elliptic counterpart, which shares the elliptic-perturbation viewpoint of \cite{JinLuQuanZhou:2025}. Note that the elliptic counterpart has been thoroughly investigated in terms of mathematical analysis \cite{Richter:1981,Alessandiri:1986,Bonito:2017} and numerical analysis \cite{Falk:1983,JinZhou:2021SINUM,Boulakia:2026}.

In the context of the subdiffusion model \eqref{eqn:fde}, early contributions on diffusion coefficient identification  mainly focus on lateral observations; see, e.g., \cite{cheng2009uniqueness,KianLiLiu:2021} for uniqueness results. The recovery from a single terminal observation poses distinct challenges and remains comparatively less understood. The only work on the terminal data is \cite{JinLuQuanZhou:2025}. Jin et al \cite{JinLuQuanZhou:2025} established a conditional stability estimate for a nontrivial spatially dependent source and derived error estimates for the recovered diffusion coefficient based on the Galerkin finite element method in space and convolution quadrature in time.
Their stability analysis essentially exploits the convergence of the terminal state $u(q)(T)$ to a nontrivial steady-state solution $v_q$, which satisfies
\begin{equation}\label{eqn:steady-state}
    \left\{\begin{aligned}
    -\nabla\cdot(q\nabla v_q)&=f
    \quad\text{in }\Omega,\\
     v_q&=0
    \quad\text{on }\partial\Omega.
    \end{aligned}\right.
\end{equation}
The time-independence of the source $f$ plays a crucial role in deriving the steady-state, and
moreover, its positivity and time-independence ensure the non-degeneracy of the limiting elliptic problem, which plays an essential role in the stability analysis via an energy argument.
The terminal-data inverse problem can then be treated as a perturbation of the elliptic problem. The key perturbation estimate takes the following form
\begin{equation*}
\big\|\partial_t^\alpha
    (u(q_1)-u(q_2))(T)\big\|_{L^2(\Omega)}
    \le c(T)\|q_1-q_2\|_{L^2(\Omega)},
    \quad \text{with}~~ c(T)\to0 \quad\text{as }T\to\infty.
\end{equation*}
This estimate allows absorbing the time-dependent contribution into the coercive elliptic estimate for sufficiently large $T$, and then using the analysis techniques for the elliptic counterpart \cite{Bonito:2017,JinZhou:2021SINUM}.

For the source-free subdiffusion model \eqref{eqn:fde}, with a zero Dirichlet boundary condition, these existing analysis strategies \cite{Triki:2021,JinLuQuanZhou:2025} do not extend directly. 
First, the steady-state argument in \cite{JinLuQuanZhou:2025} is no longer available: When $f\equiv0$, the steady state $v_q$ defined in \eqref{eqn:steady-state} vanishes identically and hence does not provide any information about the diffusion coefficient $q$. Hence, the decay of the perturbation alone is not sufficient: it must be controlled relative to the simultaneously degenerating leading term. Second,
the separation of modal decay rates used in \cite{Triki:2021} is also no longer available. Indeed, for any $0<\alpha<1$ and each fixed $j$, with  the Mittag--Leffler function $E_{\alpha,1}$ \cite[Theorem 3.2]{Jin:book2021}, the temporal factors satisfy 
\begin{equation*}
%\label{eqn:splitting}
E_{\alpha,1}(-\lambda_jT^\alpha)
\sim \frac{T^{-\alpha}}{\lambda_j\Gamma(1-\alpha)},
\quad\text{as }T\to\infty.
\end{equation*}
That is, all nonzero modal contributions have the same algebraic decay $T^{-\alpha}$, and increasing $T$ does not isolate the first eigenmode. Thus Triki's first-eigenmode reduction strategy \cite{Triki:2021} also does not extend to the subdiffusion setting. In sum, the present source free setting requires developing new analysis strategy.

\subsection{Main contributions}
We shall develop a new analysis strategy to tackle the inverse diffusivity problem for the model \eqref{eqn:fde}, based on the terminal data $u(q)(T)$. We exploit essentially the following important observation: the rescaled leading term is determined by an elliptic problem driven by the initial datum $u_0$:
\[
    T^\alpha u(q)(T)
    \rightarrow \frac{1}{\Gamma(1-\alpha)}w_q\quad\text{as}~~ T \to\infty,
    \quad \text{with }\left\{\begin{aligned} 
    -\nabla\cdot(q\nabla w_q)&=u_0
    \quad\text{in }\Omega, \\
    w_q&=0
    \quad\text{on }\partial\Omega.
    \end{aligned}\right.
\]
That is, the initial datum $u_0$ acts as an effective  source in the elliptic problem governing the asymptotic of the the rescaled limit $T^\alpha u(q)(T)$, which contains information about the diffusion coefficient $q$. This observation is the starting point for the stability analysis below. We shall establish conditional stability estimates and develop rigorous reconstruction schemes. Specifically, the main contributions of the work are four folded, which we describe as follows.

The first contribution is a conditional H\"{o}lder stability estimate in Theorem \ref{thm:conditional-stability} of recovering the diffusion coefficient $q$ from the terminal data $u(q)(T)$. The result is obtained by combining the operator splitting in \cite{LiYangZhou:2024} with the weighted energy framework of \cite{JinLuQuanZhou:2025}. The key step is to control the difference of the splitting remainders for two coefficients in terms of their $L^2(\Omega)$ distance. Using a resolvent perturbation argument, we show that the remainder can be absorbed into the coercivity estimate associated with the leading elliptic term for large terminal times. Specifically, under suitable regularity assumptions, we prove that there exists a threshold $T_*>0$ such that
\[
    \|q_1-q_2\|_{L^2(\Omega)}
    \le
    cT^{\frac\alpha2}
    \bigl(1+\|\nabla q_2\|_{L^2(\Omega)}\bigr)^{\frac12}
    \|\nabla(u(q_1)(T)-u(q_2)(T))\|_{L^2(\Omega)}^{\frac12},
\]
for all $T\ge T_*$, $q_1\in\mathcal A\cap W^{1,\infty}(\Omega)$, and $q_2\in\mathcal A$. The constants $c$ and $T_*$ depend only on an a priori bound on $\|q_1\|_{W^{1,\infty}(\Omega)}$, but are independent of $\|\nabla q_2\|_{L^2(\Omega)}$.
This estimate allows asymmetric coefficient regularity: the reference coefficient $q_1$ belongs to $W^{1,\infty}(\Omega)$, whereas the comparison coefficient $q_2$ requires only $H^1(\Omega)$ regularity. Moreover, unlike the  existing results for classical diffusion
\cite{Triki:2021,LinLin:2026}, Theorem \ref{thm:conditional-stability} does not require \textsl{a priori} knowledge
of the diffusion coefficient $q$ on the boundary $\partial\Omega$. The independence of the large-time threshold $T_*$ from the gradient of the comparison coefficient $q_2$ is essential for the error analysis of the fully discrete reconstruction scheme. 

Second, inspired by the stability analysis, we develop and analyze a fully discrete regularized least-squares reconstruction scheme with an $H^1$-seminorm penalty, using continuous piecewise linear finite elements in space and backward Euler convolution quadrature in time \cite{JinZhou:2023book}. However, the continuous stability estimate in Theorem \ref{thm:conditional-stability} does not apply to the error analysis of the discrete scheme directly, due to a lack of sufficient regularity of the data $z^\delta$. Instead, we adapt the stability argument to the discrete setting. This preserves the independence of the large-time threshold $T_*$ from the discretization and regularization parameters. 
In Theorem \ref{thm:weighted-reconstruction-error} and Corollary \ref{cor:coefficient-reconstruction-error}, we prove that any discrete reconstruction $q_h^*$ satisfies for any fixed terminal time $T\geq T_*$,
\begin{equation}\label{eqn:error-bound-recovery}
    \|q_h^*-q^\dagger\|_{L^2(\Omega)}
    \le c_T \left( h\gamma^{-1}\eta^2 +\min\{1,h^{-1}\eta\}\gamma^{-\frac12}\eta +\tau \right)^{\frac12}, \quad \text{with}~~
    \eta:=\delta+h^2+\tau+\gamma^{\frac12},
\end{equation}
The estimate \eqref{eqn:error-bound-recovery} incorporates the noise level $\delta$, the spatial mesh size $h$, the time step $\tau$, and the regularization parameter $\gamma$, which quantifies the combined effects of measurement noise, spatial and temporal discretization, and regularization. The threshold $T_*$ is independent of $h$, $\tau$, $\delta$, and $\gamma$, but $c_T$ may depend on $T$ but is independent of these parameters. The estimate \eqref{eqn:error-bound-recovery} provides a useful guideline for choosing the algorithmic parameters: the choice $\gamma\sim\delta^2$, $h\sim\delta^{\frac12}$, and $\tau=O(\delta)$ yields $\|q_h^*-q^\dagger\|_{L^2(\Omega)}=O(\delta^{\frac14})$ for any fixed terminal time $T\geq T_*$.

Third, we extend the analysis to the joint identification of the diffusion coefficient $q$ and the fractional order $\alpha$. This setting involves distinct new challenges due to the strong coupling between the two parameters and therefore requires additional conditions: additional regularity assumptions on the initial datum $u_0$, and the agreement of the coefficient $q$ with a known coefficient $q_*$ near a boundary portion $\Gamma$. Then we derive the following boundary identity
\begin{equation}\label{eqn:bdry-id}
    \operatorname{Tr}_{\Gamma}\bigl(A(q_*)u(q,\alpha)(T)\bigr)
    =
    \frac{T^{-\alpha}}{\Gamma(1-\alpha)}
    \operatorname{Tr}_{\Gamma}u_0.
\end{equation}
Moreover, the splitting remainder enjoys additional spatial regularity and  the elliptic operator applied to the remainder has a zero boundary trace. Thus, we can decouple the recovery of the fractional order $\alpha$ from the unknown interior coefficient $q$ and obtain Lipschitz stability for the order $\alpha$ in terms of the boundary observation $\operatorname{Tr}_{\Gamma}(A(q_*)u(q,\alpha)(T))$, cf. Proposition \ref{prop:order-stability}. By combining this estimate with the fixed-order stability for the coefficient $q$ that is uniform for $\alpha$ in compact subintervals of $(0,1)$, we establish joint uniqueness and conditional stability in Theorem \ref{thm:joint-q-stability} and Collary \ref{cor:joint-terminal-q-stability}.

Finally and fourth, we develop a two-stage fully discrete reconstruction based on the identity \eqref{eqn:bdry-id} and quantify error propagation between the two stages. Specifically, we first regularize boundary differentiation in \eqref{eqn:bdry-id} using local polynomial projection and recover the fractional order $\alpha$ by solving a scalar monotone equation; See Algorithm \ref{alg:order} for details. Under local $H^s$ regularity of the terminal state $u(q,\alpha)(T)$ and using polynomials of degree $p\ge2$, we prove that, for $r=\min\{s,p+1\}>\frac52$, the recovered order $\widehat{\alpha}_{\delta,T}$ satisfies
\begin{equation*}
    |\widehat\alpha_{\delta,H}-\alpha^\dagger| \le   c\left(MH^{r-\frac52}+\delta H^{-\frac52}\right),
\end{equation*}
where $H$ is mesh size of the local projection and $M$ is an upper bound on the local $H^s$ norm of the terminal state $u(q,\alpha)(T)$. Then using the recovered order $\widehat\alpha_{\delta,H}$, we reconstruct a fully discrete coefficient using the regularized least-squares formulation, and derive error estimates for the recovered fractional order $\widehat{\alpha}_{\delta,H}$ and the coefficient $q_h^*$ in Theorem \ref{thm:joint-discrete-coefficient-error}. % that account for both measurement noise and uncertainty in the fractional order.

The rest  of the paper is organized as follows.
In Section~\ref{sec:stability}, we establish resolvent and coefficient perturbation estimates and prove conditional stability for coefficient recovery with a known fractional order. In Section~\ref{sec:fem}, we present the fully discrete regularized least-squares method and derive weighted and $L^2(\Omega)$ error estimates for the recovered coefficient.
In Section~\ref{sec:joint-stability}, using an exact boundary trace identity, we establish stability for the fractional order and prove uniqueness and conditional stability for the joint recovery of the coefficient and the order. Finally, in Section~\ref{sec:joint-fem}, we develop a two-stage reconstruction procedure, derives error estimates for the recovered fractional order and coefficient. Throughout, $c$ denotes a generic constant that may change from one line to another. The notation $(\cdot,\cdot)$ denotes the $L^2(\Omega)$ inner product, and $\langle\cdot,\cdot\rangle$ denotes the duality pairing.

\section{Conditional stability for diffusion coefficient recovery}\label{sec:stability}

In this section, we establish conditional stability for recovering the diffusion coefficient $q$ from the terminal observation $u(q)(T)$, with the fractional order $\alpha\in(0,1)$ assumed known.

\subsection{Preliminary estimates}
For any $q\in\mathcal A$, let $A(q)$ be the positive self-adjoint operator on $L^2(\Omega)$ associated with the sesquilinear form
\begin{equation*}
    \mathfrak a_q(v,\varphi)
    :=( q\nabla v,\cdot\nabla\overline{\varphi}),
    \quad v,\varphi\in H_0^1(\Omega),
\end{equation*}
with its domain 
$$\operatorname{Dom}(A(q))
:=\left\{v\in H_0^1(\Omega): -\nabla\cdot(q\nabla v)\in L^2(\Omega)\right\},$$
with the divergence understood in the distributional sense. If $q\in W^{1,\infty}(\Omega)$, the standard elliptic regularity theory implies $\operatorname{Dom}(A(q))=H^2(\Omega)\cap H_0^1(\Omega)$. We also denote by $A(q)$ the associated bounded operator from $H_0^1(\Omega)$ to $H^{-1}(\Omega)$.

For $u_0\in L^2(\Omega)$, the solution $u(q)$ to problem \eqref{eqn:fde} admits the representation \cite[Section 6.2.1]{Jin:book2021} 
\begin{equation}\label{eqn:sol}
    u(q)(t)=F(t;q)u_0
    :=\frac{1}{2\pi i}\int_{\Gamma_{\theta,\sigma}}
    e^{zt}z^{\alpha-1}(z^\alpha+A(q))^{-1}u_0\,\mathrm{d}z,
    \quad t>0,
\end{equation}
where $\sigma>0$ and  the contour $\Gamma_{\theta,\sigma}\subset\mathbb{C}$ is given by
$\Gamma_{\theta,\sigma} = \{z\in \mathbb{C}:|z|=\sigma,|\arg(z)|\leq\theta\}\cup\{z\in \mathbb{C}:z=\rho e^{\pm i\theta},\rho\geq\sigma\}$, oriented with an increasing imaginary part.
Throughout, we fix $\theta\in(\frac{\pi}{2},\pi)$ so that $z^\alpha\in\Sigma_{\alpha\theta}\subset\Sigma_{\theta}$ for all $z\in\Sigma_{\theta}:=\{ z\in\mathbb{C}\setminus\{0\}:|\arg(z)|\leq\theta\}$. 

Our analysis relies crucially on the operator splitting in \cite{LiYangZhou:2024}. Specifically, the resolvent identity
\[
    (z^\alpha+A(q))^{-1}
    =A(q)^{-1}
    -z^\alpha(z^\alpha+A(q))^{-1}A(q)^{-1}
\]
and the solution representation \eqref{eqn:sol} gives the decomposition 
\begin{equation}\label{eqn:splitting}
    u(q)(t)
    =\frac{t^{-\alpha}}{\Gamma(1-\alpha)}w_q
    +R(t;q)u_0,
    \quad \text{with}\quad w_q:=A(q)^{-1}u_0,
\end{equation}
where the remainder operator $R(t;q)$ is defined by
\begin{equation}\label{eqn:remainder}
    R(t;q)
    :=-\frac{1}{2\pi i}
    \int_{\Gamma_{\theta,\sigma}}
    e^{zt}z^{2\alpha-1}
    (z^\alpha+A(q))^{-1}A(q)^{-1}\,\mathrm{d}z.
\end{equation}
The decomposition \eqref{eqn:splitting} isolates an elliptic profile $w_q$, which by definition satisfies
\[
    \left\{\begin{aligned}
    -\nabla\cdot(q\nabla w_q)&=u_0
    \quad\text{in }\Omega,\\
    \quad w_q&=0
    \quad\text{on }\partial\Omega.
    \end{aligned}\right.
\]
Thus, the rescaled terminal state $T^\alpha u(q)(T)$ takes the form
\begin{equation}\label{eqn:rescale-repres}
    T^\alpha u(q)(T)
    =\frac{1}{\Gamma(1-\alpha)}w_q
    +T^\alpha R(T;q)u_0.
\end{equation}
The splitting \eqref{eqn:rescale-repres} forms the basis of the stability analysis below: using the elliptic profile $w_q$, we establish the nondegeneracy of a rescaled weight and absorb the remainder $R(T;q)u_0$ into the weighted energy estimate for large $T$.

Now we recall two elliptic regularity estimates for $A(q)^{-1}$ \cite[Lemma 2.1]{LiSun:2017}.
%The $W^{1,\infty}(\Omega)$ estimate provides the gradient bound needed to control coefficient perturbations in $L^2(\Omega)$. 
\begin{lemma}\label{lem:elliptic} 
If $q\in W^{1,\infty}(\Omega)\cap\mathcal A$,
then there exists $c=c(\Omega, d, c_0,c_1,\|q\|_{W^{1,\infty}(\Omega)})$ such that
\begin{align*}
\|A(q)^{-1}f\|_{H^2(\Omega)} &\le c\|f\|_{L^2(\Omega)},\quad \forall f\in L^2(\Omega),\\	\|A(q)^{-1}f\|_{W^{1,\infty}(\Omega)}&\le c\|f\|_{L^\infty(\Omega)},\quad \forall f\in L^\infty(\Omega).
\end{align*}
\end{lemma}

The next result collects two useful estimates on the resolvent operator
$$
S_q(z):=(z+A(q))^{-1},\quad z\in \Sigma_{\theta}.
$$
The energy estimates follow from the spectral theorem and
duality, while the maximum-norm estimate follows from
\cite[Theorems~1.2 and~3.1(1)]{ouhabaz1995gaussian}.

\begin{lemma}\label{lem:maximum-resolv}
Fix $0<\theta<\pi$. Then for every $q\in\mathcal A$ and
$z\in\Sigma_\theta$, there exists $c=c(\Omega,d,\theta,c_0,c_1)$ such that 
\begin{align}\label{eqn:L2-resolvents}
    &\|S_q(z)\|_{L^2(\Omega)\to H_0^1(\Omega)}
    +\|S_q(z)\|_{H^{-1}(\Omega)\to L^2(\Omega)}
    \le c\min\{1,|z|^{-1/2}\},&\\
\label{eqn:max-resolvents}
&\|S_q(z)\|_{L^\infty(\Omega)\to L^\infty(\Omega)}
    \le c\min\{1,|z|^{-1}\}.&
\end{align}
%The constant $C$ depends only on $\Omega$, $d$, $\theta$, and the uniform ellipticity bounds, and is independent of $q$ and $z$.
\end{lemma}

Next we derive an $L^2(\Omega)$ perturbation estimate for the resolvent $S_q(z)$ under the asymmetric regularity assumptions on $q_1$ and $q_2$. The argument combines the resolvent bounds in Lemma \ref{lem:maximum-resolv} with elliptic regularity and requires additional regularity only for $q_1$.
\begin{lemma}
\label{lem:resdiff}
For any $q_1\in W^{1,\infty}(\Omega)\cap\mathcal A$, $q_2\in\mathcal A$, $v\in L^\infty(\Omega)$, and $z\in\Sigma_\theta$, there exists $c=c(\Omega,\theta,c_0,c_1,\|q_1\|_{W^{1,\infty}(\Omega)})$ such that
\begin{equation*} 
    \|(S_{q_1}(z)-S_{q_2}(z))v\|_{L^2(\Omega)}
    \le c|z|^{-\frac12}\|v\|_{L^\infty(\Omega)} \|q_1-q_2\|_{L^2(\Omega)}.
\end{equation*}
\end{lemma}

\begin{proof}
Let $u=S_{q_1}(z)v$.
The identity $A(q_1)u=v-zu$, Lemma~\ref{lem:elliptic} and \eqref{eqn:max-resolvents} yield
\begin{equation*} 
    \|\nabla u\|_{L^\infty(\Omega)} \le c\|v-zu\|_{L^\infty(\Omega)}\le c\|v\|_{L^\infty(\Omega)}.
\end{equation*}
Now for any $\varphi\in H_0^1(\Omega)$, we have
\begin{align*}
     \big|\langle (A(q_2)-A(q_1))u,\varphi\rangle\big|
    =\left|((q_2-q_1)\nabla u, \nabla{\varphi})      \right|
    \le \|q_1-q_2\|_{L^2(\Omega)}
          \|\nabla u\|_{L^\infty(\Omega)}
          \|\nabla\varphi\|_{L^2(\Omega)}.
\end{align*}
Combining the preceding two estimates gives 
\begin{equation}\label{eqn:A_2-A_1}
    \|(A(q_2)-A(q_1))u\|_{H^{-1}(\Omega)}
    \le c\|q_1-q_2\|_{L^2(\Omega)}
          \|v\|_{L^\infty(\Omega)}.
\end{equation}
Then the identity $ (S_{q_1}(z)-S_{q_2}(z))v
=S_{q_2}(z)(A(q_2)-A(q_1))u$, and the estimates \eqref{eqn:A_2-A_1} and \eqref{eqn:L2-resolvents} imply
\begin{align*}
   & \|(S_{q_1}(z)-S_{q_2}(z))v\|_{L^2(\Omega)}  =  \|S_{q_2}(z)(A(q_2)-A(q_1))u\|_{L^2(\Omega)}\\
    \le& c|z|^{-\frac12}
          \|(A(q_2)-A(q_1))u\|_{H^{-1}(\Omega)} \le c|z|^{-\frac12}
          \|q_1-q_2\|_{L^2(\Omega)}
          \|v\|_{L^\infty(\Omega)},
\end{align*}
which completes the proof of the lemma.
\end{proof}

\subsection{Conditional stability estimate}
Now we establish conditional stability for the inverse conductivity problem
under the asymmetric regularity assumptions $q_1\in\mathcal A\cap W^{1,\infty}(\Omega)$ and $q_2\in\mathcal A$. This setting is motivated by the numerical analysis in Section \ref{sec:fem}, where $q_1$ represents the exact coefficient and $q_2$ a reconstructed approximation whose $H^1$-seminorm is controlled by regularization.
The dependence on $\|\nabla q_2\|_{L^2(\Omega)}$ is kept explicit in the stability estimate, while the required lower bound $T_*$ on the observation time $T$ is independent of $\|\nabla q_2\|_{L^2(\Omega)}$.

Below we denote by $V_q:=T^\alpha u(q)(T)$ the rescaled state.
With the operator $R(T;q)$ defined in \eqref{eqn:remainder}, the decomposition \eqref{eqn:splitting} gives
\begin{equation}\label{eqn:Fq}
    A(q)V_q=\frac{u_0}{\Gamma(1-\alpha)}+F_q
    \quad \text{with }    F_q:=T^\alpha A(q)R(T;q)u_0.
\end{equation}

First we bound the perturbation of the remainder source $F_q$
with respect to the diffusion coefficient $q$.
\begin{lemma}\label{lem:remainder}
Let $q_1\in W^{1,\infty}(\Omega)\cap\mathcal A$,
$q_2\in\mathcal A$, and $u_0\in L^\infty(\Omega)$.
Then for $T\ge1$, there exists $c=c(c_0,c_1,\|q_1\|_{W^{1,\infty}(\Omega)})$ such that
\begin{equation*}
    \|F_{q_1}-F_{q_2}\|_{L^2(\Omega)}
    \le cT^{-\frac{\alpha}{2}}
        \|u_0\|_{L^\infty(\Omega)}
        \|q_1-q_2\|_{L^2(\Omega)}.
\end{equation*}
%The constant $C$ depends only on the fixed problem data, the uniform ellipticity bounds, and an upper bound for $\|q_1\|_{W^{1,\infty}(\Omega)}$, and is independent of $T$, $q_1$, $q_2$, and $u_0$.
\end{lemma}
\begin{proof} 
By Lemma~\ref{lem:resdiff}, we have
\begin{equation}\label{eqn:Sol-op-ell}
    \|\bigl(S_{q_1}(z^\alpha)-S_{q_2}(z^\alpha)\bigr)u_0\|_{L^2(\Omega)} \le c|z|^{-\frac{\alpha}{2}} \|u_0\|_{L^\infty(\Omega)}\|q_1-q_2\|_{L^2(\Omega)}.
\end{equation}
Then the representation \eqref{eqn:remainder} and the estimate \eqref{eqn:Sol-op-ell} yield
\begin{align*}
    \|F_{q_1}-F_{q_2}\|_{L^2(\Omega)} & = \Big\|-\frac{T^\alpha}{2\pi i}
      \int_{\Gamma_{\theta,\sigma}}
      e^{zT}z^{2\alpha-1}
      \bigl(S_{q_1}(z^\alpha)-S_{q_2}(z^\alpha)\bigr)
      u_0\,\d z\Big\|_{L^2(\Omega)}\\
    &\le cT^\alpha
        \|u_0\|_{L^\infty(\Omega)}
        \|q_1-q_2\|_{L^2(\Omega)}
        \int_{\Gamma_{\theta,\sigma}}
        e^{T\Re(z)}
        |z|^{\frac{3}{2}\alpha-1}\,|\d z|.
\end{align*}
By taking $\sigma=T^{-1}$ in the contour $\Gamma_{\theta,\sigma}$, and using the substitution $s=rT$, we obtain 
\begin{align*}
     \int_{\Gamma_{\theta,T^{-1}}}
        e^{T\Re(z)}
        |z|^{\frac32\alpha-1}\,|\d z| &=
        2\int_{T^{-1}}^\infty
          e^{rT\cos\theta}r^{\frac32\alpha-1}\,\d r
        +T^{-\frac32\alpha}
          \int_{-\theta}^{\theta}e^{\cos\phi}\,\d\phi\\
    & =
        T^{-\frac32\alpha}
        \left(
          2\int_1^\infty
            e^{s\cos\theta}s^{\frac32\alpha-1}\,\d s
          +\int_{-\theta}^{\theta}e^{\cos\phi}\,\d\phi
        \right) \le cT^{-\frac32\alpha}.
\end{align*}
Combining the last two estimates completes the proof of the lemma.
\end{proof}

Next we establish a weighted energy estimate using the argument in \cite{Bonito:2017,JinZhou:2021SINUM,JinLuQuanZhou:2025}.
\begin{lemma}\label{lem:weight-stab}
Let $q_1\in W^{1,\infty}(\Omega)\cap\mathcal A$,
$q_2\in\mathcal A$, and $u_0\in L^\infty(\Omega)$.
Then, for $T\ge1$, there holds
\begin{align*}
    &\frac12\int_\Omega
        \left(\frac{q_1-q_2}{q_1}\right)^2
        \left[
            q_1|\nabla V_{q_1}|^2
            +\left(\frac{u_0}{\Gamma(1-\alpha)}+F_{q_1}\right)             V_{q_1} \right]\,\mathrm{d}x\\
    \le&     c\bigl(1+\|\nabla q_2\|_{L^2(\Omega)}\bigr)    \|\nabla(V_{q_1}-V_{q_2})\|_{L^2(\Omega)} +cT^{-\frac\alpha2}\|q_1-q_2\|_{L^2(\Omega)}^2,
\end{align*}
where $c$ depends only on $c_0$, $c_1$
$\|q_1\|_{W^{1,\infty}(\Omega)}$, and
$\|u_0\|_{L^\infty(\Omega)}$, and is independent of
$T$ and $\|q_2\|_{H^1(\Omega)}$.
\end{lemma}
\begin{proof}
Let 
$$h:=\frac{q_1-q_2}{q_1},\quad w:=V_{q_1}-V_{q_2}\quad\mbox{and}\quad 
g_1:=\frac{u_0}{\Gamma(1-\alpha)}+F_{q_1}.$$
Then $h\in H^1(\Omega)\cap L^\infty(\Omega)$, and there hold 
\begin{equation*}
    \|h\|_{L^\infty(\Omega)}\le c, \quad
\|h\|_{L^2(\Omega)}
\le c\|q_1-q_2\|_{L^2(\Omega)} \quad \text{and}\quad
\|\nabla h\|_{L^2(\Omega)}
\le c\bigl(1+\|\nabla q_2\|_{L^2(\Omega)}\bigr).
\end{equation*}
Next we prove a uniform bound
on $V_{q_1}$.
By the resolvent estimate \eqref{eqn:max-resolvents} and the representation \eqref{eqn:remainder} (with $\sigma=T^{-1}$), we have
\begin{align}\label{eqn:Fq1}
\|F_{q_1}\|_{L^\infty(\Omega)} \le cT^\alpha\|u_0\|_{L^\infty(\Omega)}
\int_{\Gamma_{\theta,T^{-1}}}   e^{T\Re(z)}|z|^{2\alpha-1}\,|\d z|\le cT^{-\alpha}\|u_0\|_{L^\infty(\Omega)}.
\end{align}
Then Lemma~\ref{lem:elliptic} on the equation $A(q_1)V_{q_1}=g_1$ yields
\begin{equation*}
\sup_{T\ge1}\|V_{q_1}\|_{W^{1,\infty}(\Omega)}
\le c\|u_0\|_{L^\infty (\Omega)}.
\end{equation*}
Subtracting the weak formulations for $V_{q_1}$ and $V_{q_2}$ yields
\begin{equation}\label{eqn:weighted-difference}
    ((q_1-q_2)\nabla V_{q_1},\nabla\varphi)
    =(F_{q_1}-F_{q_2},\varphi)_{L^2(\Omega)}
     - (q_2\nabla w,\nabla\varphi),\quad\forall \varphi\in H_0^1(\Omega).
\end{equation}
Let $\varphi=hV_{q_1}$. Then the preceding bounds imply
that $\varphi\in H_0^1(\Omega)$ and
\begin{align}
    \|\nabla\varphi\|_{L^2(\Omega)}&\le \|V_{q_1}\|_{L^\infty(\Omega)} \|\nabla h\|_{L^2(\Omega)} +\|h\|_{L^\infty(\Omega)} \|\nabla V_{q_1}\|_{L^2(\Omega)} \le c\bigl(1+\|\nabla q_2\|_{L^2(\Omega)}\bigr),    \label{eqn:weighted-test-gradient}\\
\|\varphi\|_{L^2(\Omega)}&\le \|V_{q_1}\|_{L^\infty(\Omega)} \|h\|_{L^2(\Omega)} \le c\|q_1-q_2\|_{L^2(\Omega)}. \label{eqn:weighted-test-l2}
\end{align}
Now by using the identity $q_1-q_2=q_1h$
and testing $A(q_1)V_{q_1}=g_1$ with
$h^2V_{q_1}\in H_0^1(\Omega)$, we obtain
\begin{align*}
    ((q_1-q_2)\nabla V_{q_1},\nabla(hV_{q_1}))
  &  =\tfrac12(h^2,q_1|\nabla V_{q_1}|^2)+\tfrac12(q_1\nabla V_{q_1},\nabla(h^2V_{q_1}))\\
    &=\tfrac12(h^2,
        q_1|\nabla V_{q_1}|^2+g_1V_{q_1}).
\end{align*}
Meanwhile, Lemma~\ref{lem:remainder} and the estimate \eqref{eqn:weighted-test-l2} imply
\begin{align*}
    |(F_{q_1}-F_{q_2},\varphi)|\le \|F_{q_1}-F_{q_2}\|_{L^2(\Omega)}
\|\varphi\|_{L^2(\Omega)}\le cT^{-\frac{\alpha}{2}}\|q_1-q_2\|_{L^2(\Omega)}^2,
\end{align*}
with $c$ independent of $\|\nabla q_2\|_{L^2(\Omega)}$.
Further, the box constraint on $q_2$ and
\eqref{eqn:weighted-test-gradient} imply
\begin{align*}
    \left|(q_2\nabla w,\nabla\varphi)\right| &\le \|q_2\|_{L^\infty(\Omega)} \|\nabla w\|_{L^2(\Omega)} \|\nabla\varphi\|_{L^2(\Omega)}\le c \bigl(1+\|\nabla q_2\|_{L^2(\Omega)}\bigr) \|\nabla w\|_{L^2(\Omega)}.
\end{align*}
Combining these estimates with
\eqref{eqn:weighted-difference} proves the assertion.
\end{proof} 

To obtain an $L^2(\Omega)$ stability estimate, we impose the following assumption.
\begin{assumption}\label{ass:u_0} $\Omega\subset\mathbb R^d$ is a bounded domain
with a $C^{1,1}$ boundary, and $u_0\in L^\infty(\Omega)$ satisfies
$0<c_2\le u_0(x)\le c_3$ for  a.e. $x\in\Omega$,
for some $c_2,c_3>0$.
\end{assumption}

Next, we show that the weight in Lemma~\ref{lem:weight-stab}
has a positive lower bound for sufficiently large $T$.

\begin{lemma}\label{lem:positivity}
Let Assumption~\ref{ass:u_0} hold,
and $q\in W^{1,\infty}(\Omega)\cap\mathcal A$.
Then there exist $T_0\ge1$ and $\kappa>0$, dependent only on $d$, $\Omega$,
$\alpha$, $c_0$, $c_1$, $c_2$, $c_3$, 
and $\|q\|_{W^{1,\infty}(\Omega)}$,
and independent of $T$, such that for every $T\ge T_0$,
\begin{equation*}
    q(x)|\nabla V_q(x)|^2
    +\left(\frac{u_0(x)}{\Gamma(1-\alpha)}+F_q(x)\right)V_q(x)
    \ge\kappa
    \quad\text{for a.e. }x\in\Omega.
\end{equation*}
\end{lemma}

\begin{proof}
Let  $g_q:=\frac{u_0}{\Gamma(1-\alpha)}+F_q $. We claim
that $g_q$ is positive for large $T$. Indeed, the estimate \eqref{eqn:Fq1} implies
\begin{equation*}
    \|F_q\|_{L^\infty(\Omega)}
    \le cT^{-\alpha}\|u_0\|_{L^\infty(\Omega)}
    \le cc_3T^{-\alpha}.
\end{equation*}
Thus, with $m:=\frac{c_2}{2\Gamma(1-\alpha)}$ and $M:=\frac{c_3}{\Gamma(1-\alpha)}+m$, we can choose $T_0\ge1$ such that for all $T\geq T_0$
\begin{equation}\label{eq:positivity-source}
    m\le g_q\le M
    \quad\text{a.e. in }\Omega.
\end{equation}
Since $V_q\in H_0^1(\Omega)$ satisfies  $A(q)V_q=g_q$, the weak maximum principle implies $V_q\ge0$. Next, we control the gradient $\nabla V_q$ near the boundary $\partial\Omega$. Fix
$p>\max\{d,2\}$, and let $\beta:=1-\frac{d}{p}$.  The elliptic regularity \cite[Theorem~9.15 and Lemma~9.17]{GilbargTrudinger:2001} and Sobolev embedding $W^{2,p}(\Omega)\hookrightarrow C^{1,\beta}(\overline{\Omega})$ yield
\begin{equation*}
    \|V_q\|_{C^{1,\beta}(\overline\Omega)}    \le c\|V_q\|_{W^{2,p}(\Omega)} \le c\|g_q\|_{L^p(\Omega)} \le c_v,
\end{equation*}
where $c_v\ge1$ is independent of $T\ge T_0$. Consequently,
\begin{equation}\label{eq:positivity-gradient-holder}
    |\nabla V_q(x)-\nabla V_q(y)| \le c_v|x-y|^\beta,  \quad\forall x,y\in\overline\Omega.
\end{equation}
Next we construct a local barrier to bound $V_q$ from below. Let $D:=\operatorname{diam}(\Omega)$, fix an upper bound
$L_*$ for $\|\nabla q\|_{L^\infty(\Omega)}$. Let $b_*:=\frac{m}{2(dc_1+DL_*)}$. For any ball $B_r(z)\subset\Omega$, the function $\psi_{z,r}(x):=b_*(r^2-|x-z|^2)$ satisfies
\begin{equation*}
    A(q)\psi_{z,r}    =2b_*\bigl(dq+\nabla q\cdot(x-z)\bigr) \le 2b_*(dc_1+DL_*)    =m\le g_q \quad\text{a.e. in }B_r(z).
\end{equation*}
Since $\psi_{z,r}=0\le V_q$ on $\partial B_r(z)$, the comparison principle gives
\begin{equation}\label{eq:positivity-ball-barrier}
    V_q(x)\ge b_*(r^2-|x-z|^2),    \quad x\in B_r(z).
\end{equation}
The $C^{1,1}$ boundary admits a uniform interior ball radius $r_*>0$. For each $y\in\partial\Omega$, let $z_y:=y-r_*\nu(y)$, where $\nu(y)$ is the outward unit normal. The estimate \eqref{eq:positivity-ball-barrier} in $B_{r_*}(z_y)\subset\Omega$ yields
\begin{equation*}
    V_q(y-t\nu(y))\ge b_*(2r_*t-t^2),    \quad 0<t<r_*.
\end{equation*}
Since $V_q(y)=0$, dividing by $t$ and letting $t\to0^+$ gives 
$$-\partial_\nu V_q(y)\ge2b_*r_*.$$ 
Hence, with $\eta=2b_*r_*$,
\begin{equation}\label{eq:positivity-boundary-gradient}
    |\nabla V_q(y)|\ge\eta,    \quad y\in\partial\Omega.
\end{equation}
Next we extend the lower bound to a boundary layer: we choose the thickness $\rho:=\min\{\frac{r_*}{2},(\frac{\eta}{2c_v})^{\frac1\beta}\}$. If $\operatorname{dist}(x,\partial\Omega)<\rho$, we choose a boundary point $y$.
From the estimates \eqref{eq:positivity-gradient-holder} and
\eqref{eq:positivity-boundary-gradient},
we deduce $|\nabla V_q(x)|\ge\frac\eta2$. 
Together with $V_q\ge0$ and \eqref{eq:positivity-source}, this implies
\begin{equation*}
    q|\nabla V_q|^2+g_qV_q
    \ge\frac{c_0\eta^2}{4}
    \quad\text{for a.e. }x\in
    \{x\in\Omega:\operatorname{dist}(x,\partial\Omega)<\rho\}.
\end{equation*}
If $x\in\Omega$ satisfies
$\operatorname{dist}(x,\partial\Omega)\ge\rho$, then
$B_{\frac\rho2}(x)\subset\Omega$.
Evaluating \eqref{eq:positivity-ball-barrier} at the center
of $B_{\frac\rho2}(x)$ gives $V_q(x)\ge b_*\frac{\rho^2}{4}$. Consequently,
\begin{equation*}
    q|\nabla V_q|^2+g_qV_q
    \ge mV_q
    \ge\frac{mb_*\rho^2}{4}
    \quad\text{for a.e. }x\in
    \{x\in\Omega:\operatorname{dist}(x,\partial\Omega)\ge\rho\}.
\end{equation*}
Combining the preceding two estimates gives the desired assertion with
$\kappa:=\min\{
\frac{c_0\eta^2}{4},
        \frac{mb_*\rho^2}{4}\}>0$.
\end{proof}

By combining 
Lemma~\ref{lem:weight-stab} with
Lemma~\ref{lem:positivity}, we obtain the following conditional stability.

\begin{theorem}\label{thm:conditional-stability}
Let Assumption~\ref{ass:u_0} hold.  Fix $M_1>0$.
Then there exist $T_*\ge1$ and $c>0$ such that,
for every $q_1\in W^{1,\infty}(\Omega)\cap\mathcal A$
with $\|q_1\|_{W^{1,\infty}(\Omega)}\le M_1$,
every $q_2\in\mathcal A$, and every $T\ge T_*$,
\begin{equation*}
    \|q_1-q_2\|_{L^2(\Omega)}
    \le
    cT^{\frac\alpha2}
    \bigl(1+\|\nabla q_2\|_{L^2(\Omega)}\bigr)^{\frac12}
    \|\nabla(u(q_1)(T)-u(q_2)(T))\|_{L^2(\Omega)}^{\frac1/}.
\end{equation*}
The constants $T_*$ and $c$ depend only on $d$, $\Omega$,
$\alpha$, $c_0$, $c_1$, $c_2$, $c_3$, and $M_1$,
and are independent of $T$ and
$\|q_2\|_{H^1(\Omega)}$.
\end{theorem}

\begin{proof}
By Lemma~\ref{lem:positivity}, there exist
$T_0\ge1$ and $\kappa>0$, uniform for
$\|q_1\|_{W^{1,\infty}(\Omega)}\le M_1$, such that
\begin{equation*}
    q_1|\nabla V_{q_1}|^2
    +\left(\frac{u_0}{\Gamma(1-\alpha)}+F_{q_1}\right)V_{q_1}
    \ge\kappa
    \quad\text{a.e. in }\Omega,
    \quad T\ge T_0.
\end{equation*}
Since $q_1\le c_1$ almost everywhere in $\Omega$, Lemma~\ref{lem:weight-stab} then gives
\begin{equation*}
    \frac{\kappa}{2c_1^2}\|q_1-q_2\|_{L^2(\Omega)}^2 \le \tilde c_1\bigl(1+\|\nabla q_2\|_{L^2(\Omega)}\bigr)\|\nabla(V_{q_1}-V_{q_2})\|_{L^2(\Omega)} + \tilde c_2T^{-\frac{\alpha}{2}}\|q_1-q_2\|_{L^2(\Omega)}^2,
\end{equation*}
where $\tilde c_1,\tilde c_2>0$ depend only on the stated data and $M_1$, and are independent of $T$ and $\|q_2\|_{H^1(\Omega)}$. Then we choose the time threshold  $T_*:=\max\{T_0,(\frac{4c_1^2\tilde c_2}{\kappa})^{\frac{2}{\alpha}}\}.$
For every $T\ge T_*$, we have $\tilde c_2T^{-\frac\alpha2}\le\frac{\kappa}{4c_1^2}$.
Then using the relation $V_{q_i}=T^\alpha u(q_i)(T)$, $i=1,2$, we obtain the desired assertion:
\begin{align*}
    \frac{\kappa}{4c_1^2}
    \|q_1-q_2\|_{L^2(\Omega)}^2
    &\le
    \tilde c_1\bigl(1+\|\nabla q_2\|_{L^2(\Omega)}\bigr)
    \|\nabla(V_{q_1}-V_{q_2})\|_{L^2(\Omega)}\\
    &=
    \tilde c_1T^\alpha
    \bigl(1+\|\nabla q_2\|_{L^2(\Omega)}\bigr)
    \|\nabla(u(q_1)(T)-u(q_2)(T))\|_{L^2(\Omega)}.
\end{align*}
The construction of $T_*$ shows that it is independent
of $\|q_2\|_{H^1(\Omega)}$.
\end{proof}

\begin{remark}\label{rem:alpha-limit}
The restriction $0<\alpha<1$ is essential to the elliptic-profile splitting \eqref{eqn:splitting}.
The remainder has additional spatial smoothing and faster algebraic decay  \cite{LiYangZhou:2024}. The leading term is an elliptic profile retaining all spectral components present in $u_0$, rather than a projection onto the lowest eigenspace. Under the assumption on $u_0$, the splitting \eqref{eqn:splitting} provides the positive elliptic source and the small remainder needed in the stability analysis. The dominance of the leading term is not uniform as $\alpha\uparrow1$. Indeed, for a fixed eigenvalue $\lambda>0$, we define the relative remainder
$\varepsilon_\alpha(\lambda,T):=\lambda\Gamma(1-\alpha)T^\alpha E_{\alpha,1}(-\lambda T^\alpha)-1$.
The asymptotic expansion of the Mittag--Leffler function $E_{\alpha,1}(-t)$ \cite[Theorem 3.2]{Jin:book2021} and the convergence $E_{\alpha,1}(-\lambda T^\alpha)\to e^{-\lambda T}$ as $\alpha\uparrow1$ give, respectively,
\begin{equation*}
    \lim_{T\to\infty}\varepsilon_\alpha(\lambda,T)=0\quad\text{for each fixed }\alpha\in(0,1),\quad \lim_{\alpha\uparrow1}\varepsilon_\alpha(\lambda,T)=+\infty\quad\text{for each fixed }T>0.
\end{equation*}
Thus, even for a single eigenmode, the relative remainder cannot remain uniformly small near $\alpha=1$ for all times beyond an $\alpha$-independent threshold.  
\end{remark}

\begin{remark}\label{rem:compare}
When $\alpha=1$, the nonzero algebraically decaying elliptic leading term no longer exists. Triki \cite{Triki:2021} employed an approach based on spectral decomposition. Let $\lambda_1(q)$ be the first Dirichlet eigenvalue of $A(q)$ and $P_1(q)$ the associated spectral projection. Then
\begin{equation*}
    u(q)(T)=e^{-\lambda_1(q)T}P_1(q)u_0+e^{-TA(q)}(I-P_1(q))u_0.
\end{equation*}
Consequently,
\begin{equation*}
    A(q)u(q)(T)=\lambda_1(q)u(q)(T)+R_{\mathrm{p}}(T;q),
\end{equation*}
with $R_{\mathrm{p}}(T;q):=(A(q)-\lambda_1(q)I)e^{-TA(q)}(I-P_1(q))u_0$.
The remainder $R_{\mathrm{p}}(T;q)$ contains only higher spectral modes, and its decay relative to the first-mode is governed by the spectral gap. The stability argument additionally requires quantitative control of the dependence of the spectral quantities.
Under certain spectral, admissibility, and initial-data assumptions, Triki \cite{Triki:2021} obtained large-time Lipschitz stability from terminal data measured in $H^2(\Omega)$, with an exponentially growing stability constant. This analysis strategy differs markedly from the smoothing-based splitting argument in this work.
\end{remark}

\section{Numerical approximation and error analysis}
\label{sec:fem}

In this section, we develop and analyze a fully discrete regularized least-squares method for recovering the diffusion coefficient $q$ from the noisy terminal observation $z^\delta$ when the fractional order $\alpha\in(0,1)$ is known. 
%The state is discretized using continuous piecewise linear finite elements in space and backward Euler convolution quadrature in time. The analysis combines the minimality of discrete minimizers with the weighted stability estimate established in Section~\ref{sec:stability}. We first estimate the terminal state approximation error associated with an admissible approximation of the exact coefficient and exploit the minimizing property to derive a priori bounds for discrete minimizers. We then control the residual in the rescaled terminal equation and establish a weighted reconstruction error estimate. Finally, under the positivity assumption on the initial data, we obtain an $L^2(\Omega)$ error estimate for sufficiently large observation times.

\subsection{A fully discrete reconstruction scheme}
\label{subsec:discrete-reconstruction}
Now we develop a fully discrete scheme for reconstructing the diffusion coefficient $q$ based on standard Tikhonov regularization \cite{engl1996regularization,ItoJin:2015}.
We first describe the spatial discretization. Let $\mathcal T_h$ be a shape-regular, quasi-uniform simplicial triangulation of $\Omega$ with mesh size $h$, and on $\mathcal T_h$, we define the finite element spaces
\begin{equation*}
	V_h:=\left\{
	v_h\in C(\overline\Omega):
	v_h|_K\in\mathbb P_1(K)
	\text{ for all }K\in\mathcal T_h
	\right\}\quad\mbox{and}
	\quad
	X_h:=V_h\cap H_0^1(\Omega),
\end{equation*}
for approximating the diffusion coefficient $q$  and the state $u(q)(T)$, respectively. For domains with curved boundaries, we employ conforming finite element spaces on boundary-fitted meshes. 
The discrete admissible set $\mathcal{A}_h$ is given by
$\mathcal A_h:=\mathcal A\cap V_h$.
We also define the $L^2(\Omega)$ orthogonal projection $P_h: L^2(\Omega)\to X_h$ by
\begin{equation*}
    (P_hv,\varphi_h)=(v,\varphi_h), \quad \forall\,\varphi_h\in X_h.
\end{equation*}
The projection operator $P_h$ satisfies the following error estimate: for any $s\in[1,2]$,
\begin{equation}\label{eq:L2-projection-error}
	\|\varphi-P_h\varphi\|_{L^2(\Omega)}
	+h\|\nabla(\varphi-P_h\varphi)\|_{L^2(\Omega)}
	\le Ch^s\|\varphi\|_{H^s(\Omega)},\quad \forall \varphi\in H_0^1(\Omega)\cap H^s(\Omega).
\end{equation}
Let $\mathcal I_h$ denote the standard nodal interpolation operator on $V_h$. For $q^\dagger\in H^2(\Omega)\cap\mathcal A$, the following error estimate holds
\begin{equation}\label{eq:coefficient-interpolation}
    \|q^\dagger-\mathcal I_hq^\dagger\|_{L^2(\Omega)}+h\|\nabla(q^\dagger-\mathcal I_hq^\dagger)\|_{L^2(\Omega)}\le ch^2\|q^\dagger\|_{H^2(\Omega)}.
\end{equation}

For the temporal discretization, for any fixed $N\ge1$, we define a uniform time grid $t_n=n\tau$, $n=0,\ldots,N$, with $\tau=T/N$. Given a sequence $\{W^n=w(t_n)\}_{n=0}^N$, we approximate the Caputo derivative $\partial_t^\alpha w$ by backward Euler convolution quadrature \cite{lubich1986discretized,JinZhou:2023book}:
\begin{equation*}
    \bar\partial_\tau^\alpha W^n
    :=\tau^{-\alpha}\sum_{j=0}^n
    b_j^{(\alpha)}(W^{n-j}-W^0),
    \quad n=1,\ldots,N,
\end{equation*}
where the convolution weights $b_j^{(\alpha)}$ are determined by the generating function
\begin{equation*}
    (1-\zeta)^\alpha
    =\sum_{j=0}^\infty b_j^{(\alpha)}\zeta^j.
\end{equation*}

For each $q_h\in\mathcal A_h$, with $U_h^0=P_hu_0$, the fully discrete states $U_h^n=U_h^n(q_h)\in X_h$, $n=1,\ldots,N$, are defined by
\begin{equation}\label{eq:discrete-state}
    \bigl(\bar\partial_\tau^\alpha U_h^n,\varphi_h\bigr)
    +(q_h\nabla U_h^n,\nabla\varphi_h)=0,
    \quad \forall\,\varphi_h\in X_h,\quad n=1,\ldots,N.
\end{equation}
%Given noisy terminal data $z^\delta$ and a regularization parameter $\gamma>0$, 
Then the discrete reconstruction scheme reads
\begin{equation}\label{eq:discrete-minimization}
    \min_{q_h\in\mathcal A_h}
    \left\{J_{\gamma,h,\tau}(q_h)
    :=\|U_h^N(q_h)-z^\delta\|_{L^2(\Omega)}^2
    +\gamma\|\nabla q_h\|_{L^2(\Omega)}^2\right\},
\end{equation} 
where $\gamma>0$ is the penalty parameter.
For any fixed $h$ and $\tau$, the discrete state $U_h^N(q_h)$ depends continuously on $q_h$, which follows successively from the finite-dimensional linear systems in \eqref{eq:discrete-state}. Hence, the functional $J_{\gamma,h,\tau}$ is continuous on $\mathcal A_h$. This and the compactness of $\mathcal A_h$ imply that problem~\eqref{eq:discrete-minimization} admits at least one global minimizer $q_h^*\in\mathcal{A}_h$. 

Below we derive  error bounds on the discrete approximation $q_h^*$. The main difficulty in the analysis lies in relating the discrepancy for the fully discrete terminal state $U_h^N(q_h)$ and the noisy data $z^\delta$ controlled by $J_{\gamma,h,\tau}$  to the stability estimate in Theorem \ref{thm:conditional-stability}. Instead of using Theorem \ref{thm:conditional-stability} directly, we resort to the weighted identity underlying its proof. 
This strategy retains the dependence on the factor $\|\nabla q_h^*\|_{L^2(\Omega)}$ explicitly. One key point is that the constant of the remainder
$T^{-\frac\alpha2}
\|q^\dagger-q_h^*\|_{L^2(\Omega)}^2$
is independent of $\|\nabla q_h^*\|_{L^2(\Omega)}$ so that it can be absorbed for sufficiently large $T$.

\subsection{A priori estimates}
\label{subsec:apriori}

The next result gives an a priori error bound on the state approximation.

\begin{lemma}\label{lem:interpolated-coefficient-state-error}
Let $u_0\in L^\infty(\Omega)$, and $q^\dagger\in H^2(\Omega)\cap W^{1,\infty}(\Omega)\cap\mathcal A$. Let $u(q)$ solve \eqref{eqn:fde} and let $U_h^n(\mathcal I_hq^\dag)$ be the solution of the scheme \eqref{eq:discrete-state} with $q_h = \mathcal I_h q^\dag$. Then, for each fixed $T=N\tau>0$, there holds
\begin{equation}\label{eq:interpolated-coefficient-state-error}
    \|U_h^N(\mathcal I_hq^\dagger)-u(q^\dagger)(T)\|_{L^2(\Omega)}\le c_T(h^2+\tau).
\end{equation}
\end{lemma}

\begin{proof} 
The approximation property \eqref{eq:coefficient-interpolation} implies that  $\mathcal I_hq^\dagger \in \mathcal A_h$ and
\begin{equation*} 
    \|\mathcal I_hq^\dagger\|_{W^{1,\infty}(\Omega)}\le c
    \quad \text{and} \quad
    \|\mathcal I_hq^\dagger-q^\dagger\|_{L^2(\Omega)}\le ch^2.
\end{equation*}
The error estimates for the fully discrete scheme with nonsmooth initial data \cite[Theorem 2.4 and (3.36)]{JinZhou:2023book} lead to 
\begin{equation*}
    \|U_h^N(\mathcal I_hq^\dagger)-u(\mathcal I_hq^\dagger)(T)\|_{L^2(\Omega)} \le  c\bigl(h^2T^{-\alpha}+\tau T^{-1}\bigr) \|u_0\|_{L^2(\Omega)}.
\end{equation*} 
Meanwhile, the solution representation \eqref{eqn:sol} and Lemma~\ref{lem:resdiff}, with $q_1=q^\dagger$ and $q_2=\mathcal I_hq^\dagger$, yield
\begin{align*}
     \|u(\mathcal I_hq^\dagger)(T)-u(q^\dagger)(T)\|_{L^2(\Omega)}   & \le c\|\mathcal I_hq^\dagger-q^\dagger\|_{L^2(\Omega)}
    \|u_0\|_{L^\infty(\Omega)} \int_{\Gamma_{\theta,T^{-1}}} e^{T\Re(z)}|z|^{\frac\alpha2-1}\,|\d z|\\
    & \le    cT^{-\frac\alpha2}   \|\mathcal I_hq^\dagger-q^\dagger\|_{L^2(\Omega)}  \|u_0\|_{L^\infty(\Omega)} \le
    ch^2T^{-\frac\alpha2}\|u_0\|_{L^\infty(\Omega)}.
\end{align*}
Combining the preceding two estimates by the triangle inequality yields the estimate
\eqref{eq:interpolated-coefficient-state-error}.
\end{proof}

The minimizing property of $q_h^*$ now yields a bound  on the terminal discrepancy and the penalty simultaneously. 
\begin{lemma}\label{lem:discrete-minimizer-bound}
Under the assumptions of Lemma~\ref{lem:interpolated-coefficient-state-error},   for each fixed $T=N\tau>0$, 
there exists $c_T>0$ independent of $h$, $\tau$, $\delta$, and $\gamma$ such that the minimizer $q_h^*$ of the scheme \eqref{eq:discrete-minimization} satisfies
\begin{equation}\label{eq:discrete-minimizer-bound}
    \|U_h^N(q_h^*)-z^\delta\|_{L^2(\Omega)}+\gamma^{\frac12}\|\nabla q_h^*\|_{L^2(\Omega)}\le c_T\eta,\quad \mbox{with }\eta:=\delta+h^2+\tau+\gamma^{\frac12}.
\end{equation}
Consequently,
%\begin{equation}\label{eq:discrete-terminal-error}
$\|U_h^N(q_h^*)-u(T;q^\dagger)\|_{L^2(\Omega)}\le c_T\eta$.
\end{lemma}

\begin{proof}
The minimizing property of $q_h^*$ gives
\begin{align*}
     \|U_h^N(q_h^*)-z^\delta\|_{L^2(\Omega)}^2+\gamma\|\nabla q_h^*\|_{L^2(\Omega)}^2 
    \le\|U_h^N(\mathcal I_hq^\dagger)-z^\delta\|_{L^2(\Omega)}^2+\gamma\|\nabla\mathcal I_hq^\dagger\|_{L^2(\Omega)}^2.
\end{align*}
By the approximation property \eqref{eq:coefficient-interpolation} of the interpolation operator $\mathcal I_h$, we have
\begin{equation*}
    \|\nabla\mathcal I_hq^\dagger\|_{L^2(\Omega)}\le\|\nabla q^\dagger\|_{L^2(\Omega)}+ch\|q^\dagger\|_{H^2(\Omega)}\le c.
\end{equation*} Lemma~\ref{lem:interpolated-coefficient-state-error} and the noise condition $\|z^\delta-z^\dag\|_{L^2(\Omega)}\leq \delta$ imply
\begin{equation*}
    \|U_h^N(\mathcal I_hq^\dagger)-z^\delta\|_{L^2(\Omega)}\le c_T(h^2+\tau)+\delta.
\end{equation*}
Combining these estimates yields the bound \eqref{eq:discrete-minimizer-bound}:
\begin{equation*}
    \|U_h^N(q_h^*)-z^\delta\|_{L^2(\Omega)}^2+\gamma\|\nabla q_h^*\|_{L^2(\Omega)}^2\le c_T\bigl((\delta+h^2+\tau)^2+\gamma\bigr)\le c_T\eta^2.
\end{equation*}
The second assertion follows from the estimate \eqref{eq:discrete-minimizer-bound} and the triangle inequality.
\end{proof}

\subsection{Error estimates for the reconstruction}
\label{subsec:parameter-error}

Next we connect the bounds in Lemma \ref{lem:discrete-minimizer-bound} to the weighted
stability identity in Lemma \ref{lem:weight-stab}. In the analysis, we use the discrete elliptic operator
$A_h(q):X_h\to X_h$ defined by
\begin{equation*}
    (A_h(q)v_h,\chi_h)
    =(q\nabla v_h,\nabla\chi_h),
    \quad \forall v_h,\chi_h\in X_h.
\end{equation*}
The next lemma estimates the projected terminal source defect $T^\alpha A(q)u(q)(T)
    -T^\alpha A_h(q_h)U_h^N(q_h)$
for a reference coefficient $q\in W^{1,\infty}(\Omega)\cap\mathcal A$ and an arbitrary $q_h\in\mathcal A_h$. Note that the estimate is independent of
$\|\nabla q_h\|_{L^2(\Omega)}$.

\begin{lemma}\label{lem:projected-source-defect}
Let $u_0\in L^\infty(\Omega)$,
$q\in W^{1,\infty}(\Omega)\cap\mathcal A$, and
$q_h\in\mathcal A_h$, and
\begin{equation*}
    g_q:=T^\alpha A(q)u(q)(T)\quad\mbox{and}
    \quad
    G_{h,\tau}(q_h):=T^\alpha A_h(q_h)U_h^N(q_h).
\end{equation*}
Then, for $T=N\tau\ge1$ and $N\ge2$,
\begin{equation}\label{eq:projected-source-defect}
    \|P_hg_q-G_{h,\tau}(q_h)\|_{L^2(\Omega)}
    \le
    cT^{-\frac\alpha2}
    \bigl(\|q-q_h\|_{L^2(\Omega)}+h\bigr)
    +c\tau T^{-1}.
\end{equation}
The constant $c$ depend on $\Omega$, $\alpha$, $c_0$, $c_1$,
$\|u_0\|_{L^\infty(\Omega)}$, and
$\|q\|_{W^{1,\infty}(\Omega)}$, but is independent of
$h$, $\tau$, $T$, and $q_h$.
\end{lemma}
\begin{proof}
For $\xi\in\Sigma_\theta$, let
$$S_q(\xi):=(\xi+A(q))^{-1}\quad \mbox{and}\quad
S_{h,q_h}(\xi):=(\xi+A_h(q_h))^{-1}.$$ 
We claim that there exists $c$ independent of $q_h$ such that
\begin{equation}\label{eq:mixed-resolvent-defect}
    \|(S_{h,q_h}(\xi)P_h-P_hS_q(\xi))u_0\|_{L^2(\Omega)}
    \le c|\xi|^{-\frac12}\bigl(\|q-q_h\|_{L^2(\Omega)}+h\bigr).
\end{equation}
Indeed, let 
$$v:=S_q(\xi)u_0,
\quad v_h:=S_{h,q_h}(\xi)P_hu_0\quad\mbox{and}\quad     \vartheta_h:=v_h-P_hv.$$
The identity $A(q)v=u_0-\xi v$, Lemmas~\ref{lem:elliptic} and \ref{lem:maximum-resolv}  imply
\begin{equation*}
\|v\|_{H^2(\Omega)}+\|\nabla v\|_{L^\infty (\Omega)} \le c\|u_0\|_{L^\infty(\Omega)}.
\end{equation*}
The resolvent equations for $S_q(\xi)$ and $S_{h,q_h}(\xi)$ and the $L^2(\Omega)$-orthogonality of $P_h$ imply
\begin{align*}
    \xi(\vartheta_h,\chi_h)+(q_h\nabla\vartheta_h,\nabla\chi_h)=
    ((q-q_h)\nabla v,\nabla\chi_h)+(q_h\nabla(v-P_hv),\nabla\chi_h),
    \quad \forall\chi_h\in X_h.
\end{align*}
The approximation property \eqref{eq:L2-projection-error} of $P_h$ give
\begin{align*}
    &\left|((q-q_h)\nabla v,\nabla\chi_h)
    +(q_h\nabla(v-P_hv),\nabla\chi_h)\right| \\
     \le&
    c\left(\|q-q_h\|_{L^2(\Omega)}\|\nabla v\|_{L^\infty(\Omega)}+h\|v\|_{H^2(\Omega)}\right)\|\nabla\chi_h\|_{L^2(\Omega)}\\
    \le&
    c\bigl(\|q-q_h\|_{L^2(\Omega)}+h\bigr)\|\nabla\chi_h\|_{L^2(\Omega)}.
\end{align*}
Now we take $\xi=|\xi|e^{\mathrm i\omega}$ with $|\omega|\le\theta$. Taking $\chi_h=\vartheta_h$ in the weak formulation, and multiplying the resulting identity by $e^{-\mathrm i\frac{\omega}{2}}$ yields
\begin{align*}
 \cos(\tfrac\omega2)\left(|\xi|\|\vartheta_h\|_{L^2(\Omega)}^2
    +(q_h\nabla\vartheta_h,\nabla\vartheta_h)\right) \le
    c\bigl(\|q-q_h\|_{L^2(\Omega)}+h\bigr)\|\nabla\vartheta_h\|_{L^2(\Omega)}.
\end{align*}
Since $\cos(\tfrac\omega2)\ge\cos(\tfrac\theta2)>0$, uniform ellipticity and Young's inequality directly imply the estimate \eqref{eq:mixed-resolvent-defect}:
\begin{equation*}
    |\xi|\|\vartheta_h\|_{L^2(\Omega)}^2+\|\nabla\vartheta_h\|_{L^2(\Omega)}^2
    \le c\bigl(\|q-q_h\|_{L^2(\Omega)}+h\bigr)^2.
\end{equation*}
Let $u_h(q_h)(t)$ be the spatially semidiscrete solution with initial datum $P_hu_0$, and let
$G_h(q_h):=T^\alpha A_h(q_h)u_h(q_h)(T)$. Then the spatial error  $P_hg_q-G_h(q_h)$ can be represented by 
\begin{equation*}
    P_hg_q-G_h(q_h)
    =
    \frac{T^\alpha}{2\pi\mathrm i}
    \int_{\Gamma_{\theta,T^{-1}}}
    e^{zT}z^{2\alpha-1}
    \bigl(S_{h,q_h}(z^\alpha)P_h-P_hS_q(z^\alpha)\bigr)u_0\,\d z.
\end{equation*}
Then the estimate \eqref{eq:mixed-resolvent-defect} with $\xi=z^\alpha$ and scaling the contour with $w=Tz$ yield
\begin{align}
    \|P_hg_q-G_h(q_h)\|_{L^2(\Omega)}
    &\le cT^\alpha\bigl(\|q-q_h\|_{L^2(\Omega)}+h\bigr)
    \int_{\Gamma_{\theta,T^{-1}}}
    e^{T\Re(z)}|z|^{\frac{3}{2}\alpha-1}\,|\d z|\nonumber\\
    &\le cT^{-\frac\alpha2}\bigl(\|q-q_h\|_{L^2(\Omega)}+h\bigr).\label{eqn:err-semi}
\end{align}
Next we bound the temporal error $G_h(q_h)-G_{h,\tau}(q_h)$. Let $F_{h,q_h}(t$ and $F_{h,\tau}^n(q_h)$ denote the semidiscrete and fully discrete solution operators, respectively. Then the following representations hold  $$u_h(T;q_h)=F_{h,q_h}(T)P_hu_0\quad\mbox{and}\quad U_h^N(q_h)=F_{h,\tau}^N(q_h)P_hu_0.$$
Using the argument in \cite[Section 3.2]{JinZhou:2023book} on the kernel $K_h(z)=z^{\alpha-1}A_h(q_h)(z^\alpha+A_h(q_h))^{-1}$, we derive
\begin{equation}\label{eq:discrete-derivative-time-error}
	\bigl\|A_h(q_h)
	\bigl(F_{h,q_h}(T)-F_{h,\tau}^N(q_h)\bigr)
	\bigr\|_{L^2(\Omega)\to L^2(\Omega)}
	\le c\tau T^{-\alpha-1},
\end{equation}
with $c$ independent of $h$, $\tau$, $T$, and $q_h$. From the estimate \eqref{eq:discrete-derivative-time-error} and the $L^2(\Omega)$-stability of $P_h$, we deduce
\begin{align*}
    \|G_h(q_h)-G_{h,\tau}(q_h)\|_{L^2(\Omega)}
    =T^\alpha
    \bigl\|A_h(q_h)
    \bigl(F_{h,q_h}(T)-F_{h,\tau}^N(q_h)\bigr)
    P_hu_0\bigr\|_{L^2(\Omega)}
    \le c\tau T^{-1}\|u_0\|_{L^2(\Omega)}.
\end{align*}
This and the estimate \eqref{eqn:err-semi} complete the proof of the lemma.
\end{proof}

Next we combine the source defect estimate in Lemma \ref{lem:projected-source-defect} with the a priori bounds in Lemma \ref{lem:discrete-minimizer-bound}. Since the weighted test function used in Section~\ref{sec:stability} does not belong to $X_h$, we employ its $L^2(\Omega)$-orthogonal projection. This leads to three error terms: test-function projection error, terminal-state error, and projected source defect.  Note that the constant $c$ of the large-time remainder $cT^{-\frac\alpha2} \|q^\dagger-q_h^*\|_{L^2(\Omega)}^2$ is independent of $\|\nabla q_h^*\|_{L^2(\Omega)}$.

\begin{theorem}
\label{thm:weighted-reconstruction-error}
Let $u_0\in L^\infty(\Omega)$ and $q^\dagger\in H^2(\Omega)\cap W^{1,\infty}(\Omega)\cap\mathcal A$, and let $q_h^*$ be a minimizer of \eqref{eq:discrete-minimization}. Let $\eta$ be defined in \eqref{eq:discrete-minimizer-bound}, $
    V^\dagger:=T^\alpha u(q^\dagger)(T)$ and $g^\dagger:=\frac{u_0}{\Gamma(1-\alpha)}+F_{q^\dagger}$.
Then, for $T=N\tau\ge1$, $N\ge2$, and $0<h\le1$, there holds
\begin{align} \label{eq:weighted-reconstruction-error}
\begin{aligned}
&\quad\int_\Omega
    \left(\frac{q^\dagger-q_h^*}{q^\dagger}\right)^2
    \left(q^\dagger|\nabla V^\dagger|^2 +g^\dagger V^\dagger\right)\,\d x\\
    &\le c_T\left( h\gamma^{-1}\eta^2 +\min\{1,h^{-1}\eta\}\gamma^{-\frac12}\eta +\tau \right) +cT^{-\frac\alpha2} \|q^\dagger-q_h^*\|_{L^2(\Omega)}^2,
\end{aligned}
\end{align}
where $c_T$ is independent of $h$, $\tau$, $\delta$, and $\gamma$, and $c$ is additionally independent of $T$ and $\|\nabla q_h^*\|_{L^2(\Omega)}$.
\end{theorem}
\begin{proof}
Motivated by Lemma \ref{lem:weight-stab}, we take the test function $\varphi:=\frac{q^\dagger-q_h^*}{q^\dagger}V^\dagger$.
Note that $V^\dagger$  satisfies  
\begin{equation}\label{eqn:reg-Vdag}
    \|V^\dagger\|_{H^2(\Omega)} +\|V^\dagger\|_{W^{1,\infty}(\Omega)} \le c.
\end{equation} 
Then the box constraint on $q^\dagger$ and $q_h^*$, the regularity of $V^\dagger$ and Lemma \ref{lem:discrete-minimizer-bound} imply that $\varphi\in H_0^1(\Omega)$ and
\begin{equation}\label{eq:discrete-weighted-test-bounds}
    \|\varphi\|_{L^2(\Omega)}
    \le c\|q^\dagger-q_h^*\|_{L^2(\Omega)} 
    \quad \text{and} \quad
    \|\nabla\varphi\|_{L^2(\Omega)}
    \le c\bigl(1+\|\nabla q_h^*\|_{L^2(\Omega)}\bigr)
    \le c_T\gamma^{-\frac12}\eta.
\end{equation}
The weighted identity in Lemma \ref{lem:weight-stab} gives
\begin{equation}\label{eq:discrete-weighted-identity}
    ((q^\dagger-q_h^*)\nabla V^\dagger,\nabla\varphi)
    =
    \frac12\int_\Omega
    \left(\frac{q^\dagger-q_h^*}{q^\dagger}\right)^2
    \left(q^\dagger|\nabla V^\dagger|^2
    +g^\dagger V^\dagger\right)\,\d x.
\end{equation}
By inserting $P_h\varphi\in X_h$ and using the continuous and discrete equations at time $T$, we obtain
\begin{align*}
    ((q^\dagger-q_h^*)\nabla V^\dagger,\nabla\varphi)
    &=
    ((q^\dagger-q_h^*)\nabla V^\dagger,
    \nabla(\varphi-P_h\varphi))+
    (q_h^*\nabla(T^\alpha U_h^N(q_h^*)-V^\dagger),
    \nabla P_h\varphi)\\
    &\quad+
    (g^\dagger-G_{h,\tau}(q_h^*),P_h\varphi)=:\mathrm I_1+\mathrm I_2+\mathrm I_3.
\end{align*}
We bound the three terms separately. The box constraint of $q^\dagger$ and $q_h^*$ and the regularity in \eqref{eqn:reg-Vdag} imply
\begin{align*}
    \|\nabla\cdot((q^\dagger-q_h^*)\nabla V^\dagger)\|_{L^2(\Omega)}
    &\le
    \|q^\dagger-q_h^*\|_{L^\infty(\Omega)}
    \|\Delta V^\dagger\|_{L^2(\Omega)}+
    \|\nabla(q^\dagger-q_h^*)\|_{L^2(\Omega)}
    \|\nabla V^\dagger\|_{L^\infty(\Omega)}\\
    &\le c\bigl(1+\|\nabla q_h^*\|_{L^2(\Omega)}\bigr).
\end{align*}
Since $\varphi-P_h\varphi\in H_0^1(\Omega)$, integration by parts, the approximation property \eqref{eq:L2-projection-error} of $P_h$, and  \eqref{eq:discrete-weighted-test-bounds} give
\begin{align*}
    |\mathrm I_1|=
    \bigl|(\nabla\cdot((q^\dagger-q_h^*)\nabla V^\dagger),
    \varphi-P_h\varphi)\bigr|\le
    ch\bigl(1+\|\nabla q_h^*\|_{L^2(\Omega)}\bigr)\|\nabla\varphi\|_{L^2(\Omega)}\le c_T h\gamma^{-1}\eta^2.
\end{align*} 
Next we bound the term $\mathrm I_2$. 
The inverse inequality, the approximation property \eqref{eq:L2-projection-error} of $P_h$, Lemma~\ref{lem:discrete-minimizer-bound}, and the regularity of $V^\dagger$ in \eqref{eqn:reg-Vdag} give
\begin{align*}
\|\nabla(T^\alpha U_h^N(q_h^*)-V^\dagger)\|_{L^2(\Omega)}&\le
    \|\nabla(T^\alpha U_h^N(q_h^*)-P_hV^\dagger)\|_{L^2(\Omega)}
    +\|\nabla(P_hV^\dagger-V^\dagger)\|_{L^2(\Omega)}\\
    &\le
    ch^{-1}\|T^\alpha U_h^N(q_h^*)-P_hV^\dagger\|_{L^2(\Omega)}
+ch\|V^\dagger\|_{H^2(\Omega)}\\
&\le c_T(h^{-1}\eta+h)
    \le c_T h^{-1}\eta,
\end{align*}
Next, uniform ellipticity and the discrete smoothing estimate imply
\begin{equation*}
\|\nabla U_h^N(q_h^*)\|_{L^2(\Omega)}  \le
c\|A_h(q_h^*)^{\frac12}U_h^N(q_h^*)\|_{L^2(\Omega)}
\le cT^{-\frac\alpha2}\|u_0\|_{L^2(\Omega)}.
\end{equation*}
Together with the regularity bounds \eqref{eqn:reg-Vdag} for $V^\dagger$, this yields
$\|\nabla(T^\alpha U_h^N(q_h^*)-V^\dagger)\|_{L^2(\Omega)}\le c_T$. Combining these two bounds yields the following gradient error estimate
\begin{equation}\label{eq:discrete-terminal-gradient-error}
    \|\nabla(T^\alpha U_h^N(q_h^*)-V^\dagger)\|_{L^2(\Omega)}
    \le c_T\min\{1,h^{-1}\eta\}.
\end{equation}
Then using the $L^\infty(\Omega)$ bound on $q_h^*$, the $H^1(\Omega)$-stability of $P_h$, and \eqref{eq:discrete-weighted-test-bounds}, we deduce 
\begin{align*}
|\mathrm I_2|\le c\|\nabla(T^\alpha U_h^N(q_h^*)-V^\dagger)\|_{L^2(\Omega)}
\|\nabla P_h\varphi\|_{L^2(\Omega)}
\le c_T\min\{1,h^{-1}\eta\}\gamma^{-\frac12}\eta.
\end{align*}
Finally, by the $L^2(\Omega)$-orthogonality of $P_h$, Lemma~\ref{lem:projected-source-defect}, \eqref{eq:discrete-weighted-test-bounds} and the bound $\|q^\dagger-q_h^*\|_{L^2(\Omega)}\le c$, we obtain
\begin{align*}
|\mathrm I_3| &= |(P_hg^\dag-G_{h,\tau}(q_h^*),P_h\varphi)| \le
    c\left[T^{-\frac\alpha2}
    \bigl(\|q^\dag-q_h^*\|_{L^2(\Omega)}+h\bigr)
        +\tau T^{-1}
    \right]
    \|q^\dag-q_h^*\|_{L^2(\Omega)} \\
&\le cT^{-\frac\alpha2}\|q^\dagger-q_h^*\|_{L^2(\Omega)}^2
    +c_T(h+\tau),
\end{align*}
where the constants are independent of $T$ and $\|\nabla q_h^*\|_{L^2(\Omega)}$. Combining the preceding bounds completes the proof of the theorem.
\end{proof}

By combining Theorem \ref{thm:weighted-reconstruction-error} with Lemma~\ref{lem:positivity}, we obtain an $L^2(\Omega)$ error bound. The constant of the remainder is independent of $h$ and $\gamma$, so is the observation-time threshold.

\begin{corollary}
\label{cor:coefficient-reconstruction-error}
Let the assumptions of Theorem~\ref{thm:weighted-reconstruction-error} and Assumption~\ref{ass:u_0} hold. Then there exists $T_*\ge1$, independent of $h$, $\tau$, $\delta$ and $\gamma$, such that, for each fixed $T=N\tau\ge T_*$ with $N\ge2$ and $0<h\le1$,
\begin{equation}\label{eq:coefficient-reconstruction-error}
    \|q_h^*-q^\dagger\|_{L^2(\Omega)}\le c_T \left(
    h\gamma^{-1}\eta^2
    +\min\{1,h^{-1}\eta\}\gamma^{-\frac12}\eta
    +\tau   \right)^{\frac12},
\end{equation}
where $c_T$ is independent of $h$, $\tau$, $\delta$ and $\gamma$.
In particular, if
$\gamma\sim\delta^2$,
$h\sim\delta^{\frac12}$ and $\tau=O(\delta)$,
then, for each fixed $T\ge T_*$ and all sufficiently small $\delta>0$,
$\|q_h^*-q^\dagger\|_{L^2(\Omega)}\le c_T\delta^{\frac14}$.
\end{corollary}

\begin{proof}
By Lemma~\ref{lem:positivity}, there exist $T_0\ge1$ and $\kappa>0$, independent of $T$, $h$, $\tau$, $\delta$ and $ \gamma$, such that for all $T\ge T_0$
\begin{equation*}
    q^\dagger|\nabla V^\dagger|^2
    +g^\dagger V^\dagger
    \ge\kappa
    \quad\text{a.e. in }\Omega.
\end{equation*}
Since $0<c_0\le q^\dagger\le c_1$, we deduce
\begin{equation*}
    \int_\Omega
    \left(\frac{q^\dagger-q_h^*}{q^\dagger}\right)^2
    \left(q^\dagger|\nabla V^\dagger|^2
    +g^\dagger V^\dagger\right)\,\d x
    \ge
    \frac{\kappa}{c_1^2}
    \|q^\dagger-q_h^*\|_{L^2(\Omega)}^2.
\end{equation*}
Combining this inequality with \eqref{eq:weighted-reconstruction-error} gives
\begin{align*}
\left(c_1^{-2} \kappa -cT^{-\frac\alpha2}\right)
    \|q^\dagger-q_h^*\|_{L^2(\Omega)}^2
\le c_T\left(    h\gamma^{-1}\eta^2 +\min\{1,h^{-1}\eta\}\gamma^{-\frac12}\eta +\tau\right).
\end{align*}
Since the constant $c$ is independent of $T$, $h$, $\tau$, $\delta$, $ \gamma$ and $\|\nabla q_h^*\|_{L^2(\Omega)}$, we may choose
$$T_*=\max\left\{T_0,\,
    (\frac{2Cc_1^2}{\kappa})^{\frac2\alpha}\right\},$$
independently of $h$ and $\gamma$.  The rate $O(\delta^{\frac14})$ follows directly from the choice.
\end{proof}

\section{Simultaneous recovery of the diffusion coefficient and fractional order}
\label{sec:joint-stability}

Now we extend the stability analysis to the case of an unknown fractional order $\alpha$. Throughout, the terminal time $T$ and the initial data $u_0$ are known, and the observation is the terminal state in $\Omega$. To indicate the dependence on both unknown parameters, we denote by $u(q,\alpha)$ the solution of problem \eqref{eqn:fde}:
\begin{equation}
\label{eq:joint-forward}
    \partial_t^\alpha u(q,\alpha)+A(q)u(q,\alpha)=0,\quad 
  \mbox{with }  u(q,\alpha)(0)=u_0,
\end{equation}
where $A(q)$ is the $L^2(\Omega)$ Dirichlet realization of $-\nabla\cdot(q\nabla\cdot)$ defined in Section~\ref{sec:stability}. The terminal state  $u(q,\alpha)(T)$ is given by
\begin{equation*}
    u(q,\alpha)(T)=E_{\alpha,1}(-T^\alpha A(q))u_0,
\end{equation*}
and the exact data associated with the ground truth $(q^\dagger,\alpha^\dagger)$ is denoted by $z^\dagger=u(q^\dagger,\alpha^\dagger)(T)$.

The problem involves an unknown diffusion coefficient $q$ that determines both the eigenvalues and eigenfunctions of the elliptic operator.
Under additional regularity assumptions on the initial data $u_0$ and assuming that the coefficient $q$ is known near a portion $\Gamma$ of the boundary $\partial\Omega$, we derive a crucial boundary trace identity \eqref{eq:order-boundary-identity} that determines the fractional order $\alpha$ independently of the unknown interior coefficient $q$. This strategy enables decoupling the recovery of the two unknowns. Then we establish conditional stability for the diffusion coefficient $q$ by combining the stability estimate for the order with uniform fixed-order stability and a stability estimate of the state with respect to the order $\alpha$.

\subsection{Uniqueness and stability for order recovery}
\label{subsec:boundary-order}

Fix $I=[\underline\alpha,\overline\alpha]$, with $0<\underline\alpha<\overline\alpha<1$, and 
$T>e^{-\gamma_{\mathrm E}}$ ($\gamma_{\mathrm E}$ is
the Euler--Mascheroni constant). (We impose an  additional lower bound
on $T$ for fixed-order stability below.)
To obtain a boundary identity involving only a prescribed coefficient $q_*$ near the boundary, we require that the admissible coefficient $q$ agrees with $q_*$ in a fixed boundary neighborhood. Let $\Gamma$ be a nonempty relatively open subset  $\Gamma\subset \partial\Omega$, and let $U$ be a fixed open neighborhood of $\Gamma$. For a prescribed coefficient $q_*$, we define the following admissible set
\begin{equation}
\label{Aadj}
    \overline{\mathcal A} :=\left\{q\in\mathcal A:q=q_*\ \text{a.e. in }U\cap\Omega\right\}.
\end{equation}

\begin{assumption}
\label{ass:joint-data}
Assumption~\ref{ass:u_0} holds, and additionally $u_0\in H^2(\Omega)$.
\end{assumption}

We denote the boundary trace of the initial datum $u_0$ on $\Gamma$ by
$h_\Gamma:=\operatorname{Tr}_{\Gamma}u_0$. 
Assumption~\ref{ass:joint-data} implies
\begin{equation}\label{eqn:h-gamma}
    h_\Gamma\ge c_2
    \quad\text{a.e. on }\Gamma\quad\mbox{and}    \quad \|h_\Gamma\|_{L^2(\Gamma)}\ge c_2|\Gamma|^{\frac12}>0.
\end{equation}
For any tuple $(q,\alpha)\in\overline{\mathcal A}\times I$, the splitting \eqref{eqn:splitting} gives
\begin{equation}
\label{eq:joint-splitting}
    u(q,\alpha)(T)=\Psi_T(\alpha)w_q+R(T;q,\alpha)u_0, \quad \text{with}
    \quad \Psi_T(\alpha):=\frac{T^{-\alpha}}{\Gamma(1-\alpha)}\quad \text{and}
  \quad  w_q:=A(q)^{-1}u_0,
\end{equation}
with the remainder operator $R(T;q,\alpha)$ defined in \eqref{eqn:remainder}. For any fixed $T>0$ and $\alpha\in(0,1)$, the asymptotics of the Mittag--Leffler function  \cite[Theorem 3.2]{Jin:book2021} yields
\begin{equation*}
    E_{\alpha,1}(-T^\alpha\lambda)-\lambda^{-1}\Psi_T(\alpha) =O(\lambda^{-2}),
    \quad \mbox{as }\lambda\to\infty.
\end{equation*}
Thus, $R(T;q,\alpha)u_0\in\operatorname{Dom}(A(q)^2)$, and it follows directly from the splitting \eqref{eq:joint-splitting} that
\begin{equation*}
    A(q)u(q,\alpha)(T)=\Psi_T(\alpha)u_0+A(q)R(T;q,\alpha)u_0,
\end{equation*}
with $$A(q)R(T;q,\alpha)u_0\in\operatorname{Dom}(A(q))\subset H_0^1(\Omega).$$ 
Thus, taking the trace of $A(q)u(q,\alpha)(T)$ on $\Gamma$ eliminates the remainder $A(q)R(T;q,\alpha)u_0$. Since $q=q_*$ a.e. in the subdomain $U\cap\Omega$, we obtain the following crucial boundary trace identity
\begin{equation}
\label{eq:order-boundary-identity}
    \operatorname{Tr}_{\Gamma}\bigl(A(q_*)u(q,\alpha)(T)\bigr)=\Psi_T(\alpha)h_\Gamma,
\end{equation}
where
$A(q_*)$ denotes the local differential expression $-\nabla\cdot(q_*\nabla\cdot)$ near $\Gamma$. Its action on $u(q,\alpha)(T)$ agrees locally with the $H^1(\Omega)$ function $A(q)u(q,\alpha)(T)$. Finally, $\Psi_T$ is strictly decreasing on the interval $I$:
\begin{equation}\label{eqn:order-sign}
    \Psi_T'(\alpha)=\Psi_T(\alpha)\bigl(\psi(1-\alpha)-\log T\bigr)<0,
\end{equation}
with $\psi=\frac{\Gamma'}{\Gamma}$ is the digamma function and $\psi(1-\alpha)<\psi(1)=-\gamma_{\mathrm E}<\log T$. Together with the condition $h_\Gamma\ne0$ from \eqref{eqn:h-gamma}, the identity \eqref{eq:order-boundary-identity} uniquely determines the order $\alpha$ from the terminal observation.

The boundary identity \eqref{eq:order-boundary-identity} gives a quantitative estimate for the order $\alpha$ independently of the unknown interior coefficient $q$. In view of \eqref{eqn:order-sign}, $\Psi_T'$ is continuous and strictly negative on the interval $I$, and thus we deduce 
\begin{equation*}
m_{T,I}:=\min_{\alpha\in I}|\Psi_T'(\alpha)|> 0.
\end{equation*} 
The next result gives a stability estimate for the fractional order $\alpha$.
\begin{proposition}
\label{prop:order-stability}
Under the standing assumptions, let
$z_j:=u(q_j,\alpha_j)(T)$ with $ (q_j,\alpha_j)\in\overline{\mathcal A}\times I$, $j=1,2$. Then
\begin{equation}
\label{eq:order-trace-stability}
    |\alpha_1-\alpha_2|
    \le\frac{\|\operatorname{Tr}_{\Gamma}(A(q_*)(z_1-z_2))\|_{L^2(\Gamma)}}{m_{T,I}\|h_\Gamma\|_{L^2(\Gamma)}}
    \le\frac{\|\operatorname{Tr}_{\Gamma}(A(q_*)(z_1-z_2))\|_{L^2(\Gamma)}}{m_{T,I}c_2|\Gamma|^{\frac12}}.
\end{equation}
Moreover, let $\Gamma_0\Subset\Gamma$ be a nonempty relatively open boundary patch, and let $D\subset U\cap\Omega$ be a bounded Lipschitz subdomain such that $\Gamma_0\subset\partial D\cap\partial\Omega$. Then there exists $\tilde c_0=\tilde c_0(T, I, D, \Gamma_0,c_2)$ such that
\begin{equation}
\label{eq:order-local-stability}
    |\alpha_1-\alpha_2|\le \tilde c_0\|A(q_*)(z_1-z_2)\|_{H^1(D)}.
\end{equation}
\end{proposition}

\begin{proof}
With $v:=z_1-z_2$, the identity  \eqref{eq:order-boundary-identity} directly gives
\begin{equation*}
\operatorname{Tr}_{\Gamma}(A(q_*)v)    =\bigl(\Psi_T(\alpha_1)-\Psi_T(\alpha_2)\bigr)h_\Gamma.
\end{equation*}
Since $-\Psi_T'\ge m_{T,I}>0$ on the interval $I$, the mean value theorem yields
\begin{equation*}
    |\Psi_T(\alpha_1)-\Psi_T(\alpha_2)|\ge m_{T,I}|\alpha_1-\alpha_2|.
\end{equation*}
Consequently,
\begin{equation*}
\|\operatorname{Tr}_{\Gamma}(A(q_*)v)\|_{L^2(\Gamma)}
    =|\Psi_T(\alpha_1)-\Psi_T(\alpha_2)|\|h_\Gamma\|_{L^2(\Gamma)}
    \ge m_{T,I}|\alpha_1-\alpha_2|\|h_\Gamma\|_{L^2(\Gamma)}.
\end{equation*}
This inequality and the condition $\|h_\Gamma\|_{L^2(\Gamma)}\ge c_2|\Gamma|^{\frac12}$ from \eqref{eqn:h-gamma} show the estimate \eqref{eq:order-trace-stability}.
Next we recall that $A(q_j)z_j\in H^1(\Omega)$, $j=1,2$. Since $q_1=q_2=q_*$ a.e. in $U\cap\Omega$, we have in the sense of distribution on $D$,
\begin{equation*}
    A(q_*)v=\left.(A(q_1)z_1-A(q_2)z_2)\right|_D\in H^1(D).
\end{equation*}
%The trace of $A)q_*)v$ on $\Gamma_0$ agrees with the restriction of the boundary trace above. 
Then by restricting the boundary identity \eqref{eq:order-boundary-identity} to the part $\Gamma_0$ and the trace theorem on the subdomain $D$, we obtain
\begin{equation*}
    m_{T,I}c_2|\Gamma_0|^{\frac12}|\alpha_1-\alpha_2|
    \le\|\operatorname{Tr}_{\Gamma_0}(A(q_*)v)\|_{L^2(\Gamma_0)}
    \le c_{\mathrm{tr}}(D)\|A(q_*)v\|_{H^1(D)},
\end{equation*}
with $c_{\rm tr}(D)$ denotes the trace constant on $D$.
Thus the estimate \eqref{eq:order-local-stability} holds with $\tilde c_0=\frac{c_{\mathrm{tr}}(D)}{m_{T,I}c_2|\Gamma_0|^{\frac12}}$.
\end{proof}

\subsection{Conditional stability of the diffusion coefficient}
\label{ssec:stability2}

Next we establish conditional stability for the diffusion coefficient $q$ in the presence of an order error. The argument combines a fixed-order estimate, uniform for the order $\alpha\in I$, with a perturbation estimate for the terminal state with respect to the order $\alpha$. This step does not require the coefficients to be prescribed near the boundary patch $\Gamma$.

\begin{theorem}
\label{thm:joint-q-stability}
Let Assumption~\ref{ass:u_0} hold, and fix $M_1>0$. Let $q_1\in W^{1,\infty}(\Omega)\cap\mathcal A$, $q_2\in\mathcal A$, and $\alpha_1,\alpha_2\in I$, with $\|q_1\|_{W^{1,\infty}(\Omega)}\le M_1$. For $z_j:=u(q_j,\alpha_j)(T)$, $j=1,2$, there exist $T_*\ge1$ and $c>0$, depending only on $d$, $\Omega$, $I$, $c_0$, $c_1$, $c_2$, $c_3$, and $M_1$ but independent of $T$ and $\|\nabla q_2\|_{L^2(\Omega)}$, such that
\begin{equation}
\label{eq:joint-q-order-error}
\begin{aligned}
\|q_1-q_2\|_{L^2(\Omega)}^2
    \le c\bigl(1+\|\nabla q_2\|_{L^2(\Omega)}\bigr)\Bigl[T^{\alpha_2}\|\nabla(z_1-z_2)\|_{L^2(\Omega)}
    +(1+\log T)T^{\alpha_2-\underline\alpha}|\alpha_1-\alpha_2|\Bigr],
\end{aligned}
\end{equation}
for every $T\ge T_*$. In particular, for each fixed $T\ge T_*$,
\begin{equation}
\label{eq:joint-q-fixed-time}
\|q_1-q_2\|_{L^2(\Omega)}^2 \le c_T\bigl(1+\|\nabla q_2\|_{L^2(\Omega)}\bigr) \left(\|\nabla(z_1-z_2)\|_{L^2(\Omega)}+|\alpha_1-\alpha_2|\right),
\end{equation}
where $c_T$ may depend on $T$ but is also independent of $\|\nabla q_2\|_{L^2(\Omega)}$.
\end{theorem}
\begin{proof}
To simplify notation, we define
$$V_{q,\alpha}:=T^\alpha u(q,\alpha)(T),\quad F_{q,\alpha}:=T^\alpha A(q)R(T;q,\alpha)u_0\quad\mbox{and}\quad g_{q,\alpha}:=\frac{u_0}{\Gamma(1-\alpha)}+F_{q,\alpha}.$$
Then we have $A(q)V_{q,\alpha}=g_{q,\alpha}$. Note that the estimates in Lemma~\ref{lem:remainder} and \eqref{eqn:Fq1} hold uniformly for $\alpha\in I$. Thus, for $T\ge1$,
\begin{equation}\label{eqn:remainder-bdd}
    \|F_{q_1,\alpha}-F_{q_2,\alpha}\|_{L^2(\Omega)}
    \le cT^{-\frac{\underline\alpha}{2}}\|q_1-q_2\|_{L^2(\Omega)}\quad \text{and}
    \quad
    \|F_{q_1,\alpha}\|_{L^\infty(\Omega)}\le cT^{-\underline\alpha}
\end{equation}
with constants independent of $\|\nabla q_2\|_{L^2(\Omega)}$.
The positivity estimate in Lemma~\ref{lem:positivity} is also uniform for  $\alpha\in I$. Indeed, let
\begin{equation*}
    a_I:=\min_{\alpha\in I}\frac{1}{\Gamma(1-\alpha)}>0\quad\mbox{and}
    \quad
    b_I:=\max_{\alpha\in I}\frac{1}{\Gamma(1-\alpha)}<\infty.
\end{equation*}
The remainder bound \eqref{eqn:remainder-bdd} allows choosing $T_0\ge1$, independently of $\alpha$, such that 
for all $\alpha\in I$ and $T\ge T_0$,
\begin{equation*}
    \frac{c_2a_I}{2}\le g_{q_1,\alpha}\le c_3b_I+\frac{c_2a_I}{2}
    \quad\text{a.e. in }\Omega.
\end{equation*} 
Then the elliptic estimates and the barrier construction in Lemma~\ref{lem:positivity} provide $\kappa>0$ such that
\begin{equation*}
    q_1|\nabla V_{q_1,\alpha}|^2+g_{q_1,\alpha}V_{q_1,\alpha}\ge\kappa
    \quad\text{a.e. in }\Omega,
\end{equation*}
uniformly for $\alpha\in I$ and $T\ge T_0$.
By  Lemma~\ref{lem:weight-stab} and the preceding estimates, we obtain
\begin{equation*}
    \frac{\kappa}{2c_1^2}\|q_1-q_2\|_{L^2(\Omega)}^2
    \le \tilde c_1\bigl(1+\|\nabla q_2\|_{L^2(\Omega)}\bigr)
    \|\nabla(V_{q_1,\alpha}-V_{q_2,\alpha})\|_{L^2(\Omega)}
    +\tilde c_2T^{-\frac{\underline\alpha}{2}}\|q_1-q_2\|_{L^2(\Omega)}^2.
\end{equation*}
$\tilde c_1$ and $\tilde c_2$ are independent of $\|\nabla q_2\|_{L^2(\Omega)}$, so is the threshold
$T_*=\max\{T_0,(\frac{4c_1^2\tilde c_2}{\kappa})^{\frac{2}{\underline\alpha}}\}$. Consequently,
\begin{equation}
\label{eq:joint-uniform-fixed-order}
\begin{aligned}
    \|q_1-q_2\|_{L^2(\Omega)}^2
    \le c\bigl(1+\|\nabla q_2\|_{L^2(\Omega)}\bigr)T^\alpha
  \|\nabla\bigl(u(q_1,\alpha)(T)-u(q_2,\alpha)(T)\bigr)\|_{L^2(\Omega)},
\end{aligned}
\end{equation}
which holds uniformly for $\alpha\in I$ and $T\ge T_*$.
Next we claim that, for $\alpha,\beta\in I$ and $T\ge1$,
\begin{equation}
\label{eq:joint-forward-order-stability}
    \|\nabla\bigl(u(q_1,\alpha)(T)-u(q_1,\beta)(T)\bigr)\|_{L^2(\Omega)}
    \le c(1+\log T)T^{-\underline\alpha}|\alpha-\beta|.
\end{equation}
Indeed, with $A=A(q_1)$, the contour representation \eqref{eqn:sol} of the solution $u(q_,\alpha)$ gives
\begin{equation*}
    u(q_1,\alpha)(T)
    =\frac{1}{2\pi\mathrm i}\int_{\Gamma_{\theta,T^{-1}}}
    e^{zT}z^{\alpha-1}(z^\alpha+A)^{-1}u_0\,\d z.
\end{equation*}
%with the contour angle $\theta\in(\pi/2,\pi)$ fixed independently of $\alpha\in I$. 
Using the principal branch of the logarithm, differentiation gives
\begin{equation*}
\begin{aligned}
    \partial_\alpha\bigl[z^{\alpha-1}(z^\alpha+A)^{-1}\bigr]
    =z^{\alpha-1}\log z
    \bigl[(z^\alpha+A)^{-1}-z^\alpha(z^\alpha+A)^{-2}\bigr]=z^{\alpha-1}\log z\,A(z^\alpha+A)^{-2}.
\end{aligned}
\end{equation*}
The uniform ellipticity and Poincar\'e inequality give a lower spectral bound $\lambda_*>0$. The resolvent estimate then yields
\begin{equation*}
    \|A^{\frac12}A(z^\alpha+A)^{-2}\|_{L^2(\Omega)\to L^2(\Omega)}
    \le c\sup_{\lambda\ge\lambda_*}\lambda^{\frac32}(\lambda+|z|^\alpha)^{-2}
    \le c\lambda_*^{-\frac12},
\end{equation*}
uniformly for $\alpha\in I$. By the equivalence of $\|A^{\frac12}v\|_{L^2(\Omega)}$ and $\|\nabla v\|_{L^2(\Omega)}$ on $H_0^1(\Omega)$, setting $\xi=Tz$ and the inequality $|\log z|=|\log\xi-\log T|\le|\log\xi|+\log T$, we obtain
\begin{equation*}
\begin{aligned}
    \|\nabla\partial_\alpha u(q_1,\alpha)(T)\|_{L^2(\Omega)}
    &\le c\|u_0\|_{L^2(\Omega)}
    \int_{\Gamma_{\theta,T^{-1}}}
    e^{T\Re(z)}|z|^{\alpha-1}|\log z|\,|\d z|\le cT^{-\alpha}(1+\log T).
\end{aligned}
\end{equation*}
Integrating with respect to the order between $\alpha$ and $\beta$, and using the inequality $T^{-\gamma}\le T^{-\underline\alpha}$ for $\gamma\in I$ give the estimate \eqref{eq:joint-forward-order-stability}:
\begin{equation*}
    \|\nabla\bigl(u(q_1,\alpha)(T)-u(q_1,\beta)(T)\bigr)\|_{L^2(\Omega)}
    \le|\alpha-\beta|\sup_{\gamma\in I}
    \|\nabla\partial_\gamma u(q_1,\gamma)(T)\|_{L^2(\Omega)}
    \le c(1+\log T)T^{-\underline\alpha}|\alpha-\beta|,
\end{equation*}
with $c$ independent of $\|\nabla q_2\|_{L^2(\Omega)}$. Let $\widehat z:=u(q_1,\alpha_2)(T)$.
%The difference $\widehat z-z_2$ compares the coefficients at the common order $\alpha_2$, whereas $\widehat z-z_1$ compares the orders at the reference coefficient $q_1$. 
By the estimate \eqref{eq:joint-uniform-fixed-order} with $\alpha=\alpha_2$, we have
\begin{equation*}
    \|q_1-q_2\|_{L^2(\Omega)}^2
    \le c\bigl(1+\|\nabla q_2\|_{L^2(\Omega)}\bigr)
    T^{\alpha_2}\|\nabla(\widehat z-z_2)\|_{L^2(\Omega)}.
\end{equation*}
The triangle inequality and the estimate \eqref{eq:joint-forward-order-stability} imply
\begin{equation*}
\begin{aligned}
    \|\nabla(\widehat z-z_2)\|_{L^2(\Omega)}
    &\le\|\nabla(z_1-z_2)\|_{L^2(\Omega)}
    +\|\nabla(\widehat z-z_1)\|_{L^2(\Omega)}\\
    &\le\|\nabla(z_1-z_2)\|_{L^2(\Omega)}
    +c(1+\log T)T^{-\underline\alpha}|\alpha_1-\alpha_2|.
\end{aligned}
\end{equation*}
Then the desired estimates \eqref{eq:joint-q-order-error} and \eqref{eq:joint-q-fixed-time} follow directly.
\end{proof}

Combining Proposition~\ref{prop:order-stability} with Theorem~\ref{thm:joint-q-stability} gives the following conditional stability estimate.
\begin{corollary}
\label{cor:joint-terminal-q-stability}
Under the assumptions of Theorem~\ref{thm:joint-q-stability} and Assumption~\ref{ass:joint-data},
let $(q_j,\alpha_j)\in\overline{\mathcal A}\times I$, $j=1,2$, and fix $T\ge T_*$.  Then there holds
\begin{equation*}
    \|q_1-q_2\|_{L^2(\Omega)}
    \le c_T\bigl(1+\|\nabla q_2\|_{L^2(\Omega)}\bigr)^{\frac12}
    \left(\|\nabla(z_1-z_2)\|_{L^2(\Omega)}
    +\|\operatorname{Tr}_{\Gamma}(A(q_*)(z_1-z_2))\|_{L^2(\Gamma)}\right)^{\frac12}.
\end{equation*}
In particular, for $\Gamma_0$ and $D$ in Proposition~\ref{prop:order-stability},
\begin{equation*}
    \|q_1-q_2\|_{L^2(\Omega)}
    \le c_T\bigl(1+\|\nabla q_2\|_{L^2(\Omega)}\bigr)^{\frac12}
    \left(\|\nabla(z_1-z_2)\|_{L^2(\Omega)}
    +\|A(q_*)(z_1-z_2)\|_{H^1(D)}\right)^{\frac12}.
\end{equation*}
The constant $c_T$ may depend on $T$, but is independent of $\|\nabla q_2\|_{L^2(\Omega)}$. 
\end{corollary}

\section{Numerical approximation of joint recovery }
\label{sec:joint-fem}

The boundary trace identity \eqref{eq:order-boundary-identity} provides a natural way to construct a two-stage algorithm by decoupling the recovery of the fractional order $\alpha$ from that of the diffusion coefficient $q$. From the identity \eqref{eq:order-boundary-identity}, we first determine the order $\alpha$ via a scalar relation involving only the terminal state and the known coefficient $q_*$ near $\Gamma$, and then recover the diffusion coefficient $q$ using the regularized least-squares method in Section~\ref{sec:fem}.
Nonetheless, the development of a stable numerical procedure is challenging. First, the identity \eqref{eq:order-boundary-identity} involves second-order spatial derivatives of the terminal state, which is unstable on the noisy terminal observation $z^\delta$. To resolve the challenge, we regularize the operation by local polynomial projection near the boundary $\Gamma$ and then recover the order $\alpha$ by a scalar monotone inversion. Second, the recovered order introduces a model discrepancy into the recovery of the coefficient $q$. Fortunately, the discrepancy can be quantified via  the fixed-order error analysis.

\subsection{Numerical recovery of fractional order}
\label{ssec:joint-fem} 

First we discuss the recovery of the order $\alpha$. Throughout, we use the setting and notation of Section~\ref{sec:joint-stability}, $(q^\dagger,\alpha^\dagger)\in\overline{\mathcal A}\times I$ is the ground truth, $z^\dagger$ is the exact terminal data, and $q_*$ is the prescribed coefficient in $U\cap\Omega$, and also $\Psi_T$, $m_{T,I}$, and $h_\Gamma$ etc.

\begin{assumption}
\label{ass:geometry}
The boundary $\Gamma$ contains a flat portion. After a rigid coordinate change, let
\begin{equation*}
S=\prod_{j=1}^{d-1}(a_j,b_j),
\quad
D=S\times(0,a),
\quad
\Gamma_0=S\times\{0\},
\end{equation*}
where $a>0$ is fixed,
$D\subset U\cap\Omega$, $
\overline{\Gamma_0}\subset\Gamma$, and $\Omega$ lies locally on the side $\{x_d>0\}$.
\end{assumption}

The pair $(D,\Gamma_0)$  specifies one choice of the local observation geometry in Proposition~\ref{prop:order-stability}. The next  assumption imposes the local regularity on the exact data $z^\dag$.
%which is needed for approximate the second-order boundary quantity from noisy volume observations. The fractional order is reconstructed using only data on $D$. In particular, 

\begin{assumption}
\label{ass:data}
$q_*\in W^{2,\infty}(D)$, and for some $s>\frac52$ and $M>0$, $z^\dag$ satisfies
$z^\dagger|_D\in H^s(D)$ and $\norm{z^\dagger}{H^s(D)}\le M$.
\end{assumption}

For $v\in H^s(D)$, we can define a boundary operator
\begin{equation}
\label{eq:boundary-operator}
Bv
:=
-(\Tr_{\Gamma_0}q_*)\Tr_{\Gamma_0}(\Delta v)
-(\Tr_{\Gamma_0}\nabla q_*)\cdot
\Tr_{\Gamma_0}(\nabla v).
\end{equation}
Since $s>\frac52$, $\nabla v$ and $\Delta v$ admit $L^2(\Gamma_0)$ traces.
Moreover, Assumption~\ref{ass:data} ensures that the traces
of $q_*$ and $\nabla q_*$ belong to $L^\infty(\Gamma_0)$.
Thus, the operator $B$ is well defined.
Note that for the exact terminal state $z^\dag$, $Bz^\dag$ agrees with the restriction of the boundary trace in
\eqref{eq:order-boundary-identity}. These observations and  \eqref{eqn:h-gamma} imply
\begin{equation}\label{eq:identification}
g^\dagger:=Bz^\dagger
=
\Psi_T(\alpha^\dagger)h_\Gamma
\quad\text{in }L^2(\Gamma_0)\quad\mbox{and}\quad 
\norm{h_\Gamma}{L^2(\Gamma_0)}
\ge c_2|\Gamma_0|^{\frac12}>0.
\end{equation}
%The additional local regularity in Assumption~\ref{ass:data}, however, remains a hypothesis of the numerical error analysis.

Below we denote the range of the scalar identification map $\Psi_T$ by
\begin{equation*}
J_T:=\Psi_T(I) = [\Psi_T(\overline\alpha),\Psi_T(\underline\alpha)].
\end{equation*}
%We use an auxiliary mesh of $D$ to regularize the boundary differentiation. Its size $H$ need not coincide with the mesh size used for the forward solver or the coefficient reconstruction in Section~\ref{sec:fem}.

For the numerical recovery of the order $\alpha$, we employ an auxiliary mesh (of mesh size $H$) for a fixed observation region $D$ within $U\cap\Omega$ adjacent to a flat portion of the observation boundary $\Gamma$. 
\begin{assumption}%[Boundary-fitted affine meshes]
\label{ass:mesh}
Let $\{\tTH\}$ be a shape-regular, quasi-uniform family of conforming affine simplicial meshes of $D$. The meshes represent $D$ exactly, and $\Gamma_0$ is a union of mesh boundary faces up to a set of surface measure zero. With $
H_K:=\diam(K)$ and $H:=\max_{K\in\tTH}H_K$ $(0<H\le1$), there exists $c_{\mathrm{qu}}>0$, independent of $H$,
such that
\begin{equation*}
c_{\mathrm{qu}}H\le H_K\le H
\quad\forall K\in\tTH.
\end{equation*}
$\FH$ denotes the collection of boundary faces contained in $\overline{\Gamma_0}$. For any $F\in\FH$, $K_F$ denotes the unique adjacent mesh element.
\end{assumption}

Fix $p\ge2$, we define a discontinuous polynomial space
\begin{equation*}
W_H^p := \left\{ v_H\in L^2(D): v_H|_K\in\PP_p(K)\ \forall  K\in\tTH \right\},
\end{equation*}
where $\PP_p(K)$ denotes the space of polynomials on $K$ of total degree at most $p$. $\Pi_K:L^2(K)\to\PP_p(K)$ denotes the local $L^2(K)$-orthogonal projection defined by
\begin{equation}\label{eq:local-projection}
\ip{\Pi_Kv}{w}{L^2(K)} = \ip{v}{w}{L^2(K)},
\quad\forall w\in\PP_p(K).
\end{equation}
Each local projection is uniquely defined because its mass matrix is positive definite. 
The elementwise projection $\Pi_H:L^2(D)\to W_H^p$ is defined by $(\Pi_Hv)|_K=\Pi_K(v|_K)$. Then we have 
\begin{equation}\label{eq:projection-stability}
\norm{\Pi_Kv}{L^2(K)}\le\norm{v}{L^2(K)}\quad\mbox{and} \quad \norm{\Pi_Hv}{L^2(D)}\le\norm{v}{L^2(D)}.
\end{equation}
The projected data $v_H=P_Hv$ can be differentiated within each element $K$. Accordingly, we define a discrete analogue $B_H:W_H^p\to L^2(\Gamma_0)$ of the operator $B$ facewise by
\begin{equation}
\label{eq:discrete-boundary}
(B_Hv_H)|_F := -q_*|_F\,\bigl(\Delta(v_H|_{K_F})\bigr)|_F -(\nabla q_*)|_F\cdot \bigl(\nabla(v_H|_{K_F})\bigr)|_F, \quad F\in\FH,
\end{equation}
where the derivatives of $v_H$ are taken within $K_F$ before restriction to $F$.

The reconstruction procedure consists of four steps; see Algorithm \ref{alg:order} for details. The first two steps approximate the boundary quantity in \eqref{eq:identification}; 
the third step computes the scalar $c^\delta_H\approx \Psi_T(\alpha^\dag)$, and the last step inverts the function $\Psi_T$. For a closed interval $[a,b]\subset\R$, we define the cut-off projection $P_{[a,b]}x:=\min\{b,\max\{a,x\}\}$. The cut-off  $P_{J_T}c_H^\delta$ in Step 4 of the algorithm maps the scalar $c^\delta_H$ into $J_T$, on which the identification map $\Psi_T$ has a unique inverse. Thus the reconstruction procedure is well defined.

\begin{algorithm}
    \caption{The recovery of the order $\alpha$.\label{alg:order}}
    \begin{algorithmic}[1]
        \STATE  Solve the local projection problem 
\eqref{eq:local-projection}, and compute the projection $z_H^\delta:=\Pi_Hz^\delta$.
\STATE Compute the discrete boundary projected data $
g_H^\delta:=B_Hz_H^\delta$.

\STATE Estimate $\Psi_T(\alpha^\dagger)$ by $c_H^\delta := \norm{h_\Gamma}{L^2(\Gamma_0)}^{-2}\ip{g_H^\delta}{h_\Gamma}{L^2(\Gamma_0)}$.

\STATE Compute 
$\widehat\alpha_{\delta,H} := \Psi_T^{-1}\bigl(P_{J_T}c_H^\delta\bigr)$.

\STATE Compute the unique $\widehat\alpha_{\delta,H}\in I$.
    \end{algorithmic}
\end{algorithm}

We first recall useful local approximation estimates for the elementwise $L^2$-projection.
\begin{lemma}\label{lem:projection}
Let $\frac52<r\le p+1$, and set
\begin{equation*}
\theta:= \begin{cases}
r-2, & \frac52<r<3,\\
1, & r\ge3.
\end{cases}
\end{equation*}
Then $\frac12<\theta\le1$ and $2+\theta\le r$.
For $v\in H^r(K)$, let $e_K:=\Pi_Kv-v$. There exists $c$,
depending on $r$, $p$ and $d$, but independent of $K$ and $H$, such that
\begin{equation} \label{eqn:local-projection}
\sum_{j=0}^{2} H_K^j \|D^j e_K\|_{L^2(K)}
+ H_K^{2+\theta}|D^2 e_K|_{H^\theta(K)}
\le c H_K^r |v|_{H^r(K)}.
\end{equation}
\end{lemma}
\begin{proof}
The estimate for the $L^2$ error can be found at  \cite[Theorem~18.16(iii)]{ErnGuermond2021}. The estimate for $\theta\in (\frac12,1]$  follows from that for the third derivative when $r\ge3$, and from local quadratic approximation and the inverse inequality when $\frac52<r<3$.
\end{proof}

We also use two trace estimates. The first estimate follows from \cite[Lemma~12.15 and Remark~12.19]{ErnGuermond2021}, while the second from the discrete trace and inverse inequalities in \cite[Lemmas~12.1 and~12.10]{ErnGuermond2021}. \begin{lemma}
\label{lem:trace}
Let $K$ be a simplex  and let $F$ be a face of $K$. Then there exists $c$ depending only on $\theta$, $p$, $d$ but not on $K$ or $H$ such that
\begin{equation}\label{eq:scaled-trace}
\norm{w}{L^2(F)} \le c( H_K^{-\frac12}\norm{w}{L^2(K)}
+ H_K^{\theta-\frac12}\seminorm{w}{H^\theta(K)}),\quad \forall w\in H^\theta(K), \quad \tfrac12<\theta\le1.
\end{equation}
\begin{equation}\label{eq:polynomial-trace}
\norm{D^jw_H}{L^2(F)} \le c H_K^{-j-\frac12}\norm{w_H}{L^2(K)}, \quad \forall w_H\in\PP_p(K), \quad j=1,2. 
\end{equation}
\end{lemma}

The next result gives an error bound on the recovered order $\widehat\alpha_{\delta,H}$.
\begin{theorem} \label{thm:reconstruction}
Let Assumptions~\ref{ass:geometry}, \ref{ass:data}, and~\ref{ass:mesh} hold. Fix $p\ge2$ and let $r:=\min\{s,p+1\}>\frac52$. Then Algorithm \ref{alg:order} yields a unique $\widehat\alpha_{\delta,H}\in I$. Moreover, there exists $c$ depending on $D$, $d$, $s$, $p$, $\norm{q_*}{W^{2,\infty}(D)}$ but independent of $H$, $\delta$ and $M$ such that
\begin{align}
\label{eq:boundary-error}
\norm{g_H^\delta-g^\dagger}{L^2(\Gamma_0)}
&\le c\left(MH^{r-\frac52}+\delta H^{-\frac52}\right),\\
\label{eq:order-error}
|\widehat\alpha_{\delta,H}-\alpha^\dagger|
&\le cm_{T,I}^{-1}\norm{h_\Gamma}{L^2(\Gamma_0)}^{-1}
\left(MH^{r-\frac52}+\delta H^{-\frac52}\right).
\end{align}
\end{theorem}

\begin{proof}

Let $v_H\in W_H^p$. For each $F\in\FH$, let $K=K_F$ denote its adjacent element. The estimate \eqref{eq:polynomial-trace} implies
\[
\norm{B_Hv_H}{L^2(F)} \le c(\norm{D^2(v_H|_K)}{L^2(F)} +\norm{\nabla(v_H|_K)}{L^2(F)}) \le cH_K^{-\frac52}\norm{v_H}{L^2(K)}.
\]
Summing the squared inequality over the observation faces, and using quasi-uniformity of the mesh give
\begin{equation}
\label{eq:boundary-stability}
\norm{B_Hv_H}{L^2(\Gamma_0)}
\le c H^{-\frac52}\norm{v_H}{L^2(D)}.
\end{equation}
Then the $L^2$ stability of the projection $\Pi_H$ yields
\begin{equation}
\label{eq:noise-amplification}
\norm{B_H\Pi_H(z^\delta-z^\dagger)}{L^2(\Gamma_0)}
\le cH^{-\frac52}\norm{z^\delta-z^\dagger}{L^2(D)}
\le c\delta H^{-\frac52}.
\end{equation}
For $F\in\FH$ and $K=K_F$, let $e_K:=\Pi_K(z^\dagger|_K)-z^\dagger|_K$, and choose $\theta$ as in Lemma~\ref{lem:projection}. 
By the estimate \eqref{eq:scaled-trace} (on $D^2e_K$) and using the local projection estimates in \eqref{eqn:local-projection} gives
\begin{equation}
\label{eq:second-derivative-trace-error}
\norm{D^2e_K}{L^2(F)}\le c\left(
H_K^{-\frac12}\norm{D^2e_K}{L^2(K)}
+H_K^{\theta-\frac12}\seminorm{D^2e_K}{H^\theta(K)}
\right)\le cH_K^{r-\frac52}\seminorm{z^\dagger}{H^r(K)}.
\end{equation}
Similarly, the trace estimate in Lemma \ref{lem:trace} yields
\begin{equation} \label{eq:first-derivative-trace-error}
\norm{\nabla e_K}{L^2(F)}\le c\left( H_K^{-\frac12}\norm{\nabla e_K}{L^2(K)} +H_K^{\frac12}\norm{D^2e_K}{L^2(K)} \right)\le cH_K^{r-\frac32}\seminorm{z^\dagger}{H^r(K)}.
\end{equation}
%Since $r>5/2$, the traces of the exact derivatives are well defined and agree with the corresponding elementwise traces. 
By the identity $g^\dagger=Bz^\dagger$,and the definitions of the operators $B$ and $B_H$, we obtain
\[
\norm{B_H\Pi_Hz^\dagger-g^\dagger}{L^2(F)}
\le cH_K^{r-\frac52}\seminorm{z^\dagger}{H^r(K)}.
\]
Summing the squared inequality over the observation faces and using
$H_K\le H$ and $2r-5>0$ gives
\[
\norm{B_H\Pi_Hz^\dagger-g^\dagger}{L^2(\Gamma_0)}^2\le c H^{2r-5} \sum_{K\in\tTH}\seminorm{z^\dagger}{H^r(K)}^2\le cH^{2r-5}\norm{z^\dagger}{H^r(D)}^2 \le c M^2H^{2r-5}.
\]
Hence
\begin{equation}
\label{eq:consistency}
\norm{B_H\Pi_Hz^\dagger-g^\dagger}{L^2(\Gamma_0)}
\le C M H^{r-\frac52}.
\end{equation}
The splitting
$g_H^\delta-g^\dagger
=B_H\Pi_H(z^\delta-z^\dagger)
+\bigl(B_H\Pi_Hz^\dagger-g^\dagger\bigr)$ and the estimates \eqref{eq:noise-amplification} and \eqref{eq:consistency} imply the estimate \eqref{eq:boundary-error}. 
Next, in view of \eqref{eq:identification}, let
\[
c^\dagger
:= \norm{h_\Gamma}{L^2(\Gamma_0)}^{-2}\ip{g^\dagger}{ h_\Gamma}{{L^2(\Gamma_0)}} =\Psi_T(\alpha^\dagger)\in J_T.
\]
Then the Cauchy--Schwarz inequality gives
\begin{equation}\label{eq:scalar-error}
|c_H^\delta-c^\dagger|
\le \norm{h_\Gamma}{L^2(\Gamma_0)}^{-1}\norm{g_H^\delta-g^\dagger}{L^2(\Gamma_0)}.
\end{equation}
By Proposition~\ref{prop:order-stability}, $\Psi_T^{-1}$ is Lipschitz on $J_T$ with Lipschitz constant $m_{T,I}^{-1}$. Therefore, we have
\[
|\widehat\alpha_{\delta,H}-\alpha^\dagger| \le m_{T,I}^{-1}|c_H^\delta-c^\dagger| \le  m_{T,I}^{-1}\norm{h_\Gamma}{L^2(\Gamma_0)}^{=1} \norm{g_H^\delta-g^\dagger}{L^2(\Gamma_0)}.
\]
Substituting the estimate \eqref{eq:boundary-error} proves the assertion \eqref{eq:order-error}.
\end{proof}

In the procedure, the mesh size $H$ acts as a regularization parameter. The next result gives a sufficient condition for convergence and an a priori choice obtained by balancing approximation and noise errors. 
Piecewise quadratic polynomials provide the lowest-degree choice covered by the analysis. 

\begin{corollary}\label{cor:parameter}
Under the assumptions of Theorem~\ref{thm:reconstruction}, if the mesh size $H$ is chosen such that $ H \to0$ and $\delta H ^{-\frac52}\to0$  as $\delta\to0$, then the convergence $\widehat\alpha_{\delta,H(\delta)} \to\alpha^\dagger$ holds. Moreover, for the choice $
H\sim(\frac\delta M)^{\frac1r}$, there holds
\begin{equation}\label{eq:balanced-rate}
|\widehat\alpha_{\delta,H}-\alpha^\dagger| \le cm_{T,I}^{-1}\norm{h_\Gamma}{L^2(\Gamma_0)}^{-1} M^{\frac{5}{2r}}\delta^{1-\frac{5}{2r}}.
\end{equation}
\end{corollary}

\begin{comment}
\begin{corollary}
\label{cor:quadratic}
Suppose that the assumptions of Theorem~\ref{thm:reconstruction} hold with $p=2$ and $s\ge3$. Then $r=3$, and
\begin{equation} \label{eq:quadratic-error}
|\widehat\alpha_{\delta,H}-\alpha^\dagger| \le cm_{T,I}^{-1}\norm{h_\Gamma}{L^2(\Gamma_0)}^{-1} \left(MH^{\frac12}+\delta H^{-\frac52}\right).
\end{equation}
Choosing $H\sim(\frac\delta M)^{\frac13}$ gives
$|\widehat\alpha_{\delta,H}-\alpha^\dagger| \le cm_{T,I}^{-1}\norm{h_\Gamma}{L^2(\Gamma_0)}^{-1}
M^{\frac56}\delta^{\frac16}$.
If $\frac52<s<3$, the reconstruction remains valid with $r=s$ in Theorem~\ref{thm:reconstruction}.
\end{corollary}
\end{comment}

\subsection{Diffusion coefficient recovery with the reconstructed order}
\label{ssec:joint-coefficient-recovery}

Let $\alpha^*=\widehat\alpha_{\delta,H}\in I$ be the recovered order.
Now we recover the diffusion coefficient $q$ using the method of Section~\ref{sec:fem} with $\alpha=\alpha^*$, and employ the spaces, projection, and time grid in Section~\ref{subsec:discrete-reconstruction}. Then for $q_h\in\mathcal A_h$, with $U_h^0(q_h;\alpha^*)=P_hu_0$, we define $U_h^n(q_h;\alpha^*)\in X_h$, $n=1,\ldots,N$, by
\begin{equation}\label{eq:joint-discrete-state}
    \bigl(\bar\partial_\tau^{\alpha^*}
    U_h^n(q_h;\alpha^*),\chi_h\bigr)
    +(q_h\nabla U_h^n(q_h;\alpha^*),\nabla\chi_h)=0, \quad\forall
 \chi_h\in X_h.
\end{equation}
Here $\bar\partial_\tau^{\alpha^*}$ is the backward Euler convolution quadrature (at the order $\alpha^*$). Then the fully discrete scheme reads 
\begin{equation}
\label{eq:joint-discrete-minimization}
\min_{q_h\in\mathcal A_h}\left \{
    J_{\gamma,h,\tau}^{\alpha^*}(q_h):= \|U_h^N(q_h;\alpha^*)-z^\delta\|_{L^2(\Omega)}^2
    +\gamma\|\nabla q_h\|_{L^2(\Omega)}^2 \right\},
\end{equation}
where $\gamma>0$ is the regularization parameter.
The existence of a minimizer $q_h^*\in\mathcal{A}_h$ follows similarly as in Section~\ref{subsec:discrete-reconstruction}.  The next result gives an error bound on the approximation $q_h^*$.
%The boundary-neighborhood constraint used to identify the order is not needed for the coefficient estimate below, so we retain the admissible set $\mathcal A_h$.

\begin{theorem}\label{thm:joint-discrete-coefficient-error}
Let Assumption~\ref{ass:u_0} hold, and $q^\dagger\in H^2(\Omega)\cap W^{1,\infty}(\Omega) \cap\mathcal A$.
Let $\alpha^\dagger,\alpha^*\in I$,  and let $q_h^*$ be any  minimizer of problem \eqref{eq:joint-discrete-minimization}.
Then there exists $T_*\ge1$ such that, for each fixed $T=N\tau\ge T_*$ with $N\ge2$ and $0<h\le1$,
\begin{equation}
\label{eq:joint-discrete-coefficient-error}
    \|q_h^*-q^\dagger\|_{L^2(\Omega)}
    \le c_T
    \left(
        h\gamma^{-1}\eta_*^2
        +\min\{1,h^{-1}\eta_*\}\gamma^{-\frac12}\eta_*
        +\tau
    \right)^{\frac12},\quad \mbox{with } \eta_*:=\eta+|\alpha^\dagger-\alpha^*|.
\end{equation}
The constants $T_*$ and $c_T$ depend only on $I$, and  $q^\dagger$, and $c_T$ additionally depends on $T$.
%Both are independent of $h$, $\tau$, $\delta$, $\gamma$, $\alpha^\dagger,\alpha^*\in I$, and $\|\nabla q_h^*\|_{L^2(\Omega)}$.
\end{theorem}

\begin{proof}
The proof is based on reduction to the fixed-order analysis in Section \ref{sec:fem}, using the state $z_*^\dagger:=u(q^\dagger,\alpha^*)(T)$. By \eqref{eq:joint-forward-order-stability}, and  the Poincar\'e inequality, we have
\begin{equation*}
\begin{aligned}
    \|z^\delta-z_*^\dagger\|_{L^2(\Omega)}\le \|z^\delta-z^\dagger\|_{L^2(\Omega)}
       +\|z^\dagger-z_*^\dagger\|_{L^2(\Omega)}\le \delta
       +c(1+\log T)T^{-\underline\alpha}
        |\alpha^\dagger-\alpha^*|\le \delta+c_T|\alpha^\dagger-\alpha^*|.
\end{aligned}
\end{equation*}
Then using the argument of Lemmas~\ref{lem:interpolated-coefficient-state-error} and~\ref{lem:discrete-minimizer-bound}, we derive
\begin{equation*}
\|U_h^N(q_h^*;\alpha^*)-z_*^\dagger\|_{L^2(\Omega)}
    +\gamma^{\frac12}\|\nabla q_h^*\|_{L^2(\Omega)}
    \le c_T\eta_*,
\end{equation*}
where $c_T$ is uniform for $\alpha^*\in I$. % the contour angles can be fixed on $I$, and the time convolution estimates remain bounded because $\underline\alpha>0$ and $\overline\alpha<1$.
Thus, the argument of Theorem~\ref{thm:weighted-reconstruction-error} applies with the reference state $z_*^\dagger$ and the order $\alpha^*$. The remainder is bounded by
$ cT^{-\frac{\underline\alpha}{2}}\|q_h^*-q^\dagger\|_{L^2(\Omega)}^2 $. The positivity estimate is also uniform on the interval $I$, cf. the proof of Theorem~\ref{thm:joint-q-stability}. 
Thus, the absorption argument of Corollary~\ref{cor:coefficient-reconstruction-error} gives the desired estimate.
\end{proof}

\begin{remark}\label{rem:two-stage-error}
Under the assumptions of Theorem~\ref{thm:joint-discrete-coefficient-error}, let $e_\alpha$ be an upper bound for $|\alpha^\dagger-\alpha^*|$ and set $\rho:=\delta+e_\alpha$. For any fixed $T\ge T_*$, the choice $\gamma\sim\rho^2$, $h\sim\rho^{\frac12}$ and $\tau=O(\rho)$ 
yields $\eta_*=O(\rho)$ and 
$\|q_h^*-q^\dagger\|_{L^2(\Omega)} \le c_T\rho^{\frac14}$
as $\rho\to0$. Thus, the discretization and regularization parameters should account for both the noise level $\delta$  and the error of the order recovery.
Moreover, let $r:=\min\{s,p+1\}>\frac52$ and $\beta:=1-\frac{5}{2r}\in(0,1)$, we  
choose the mesh size $H\sim(\frac{\delta}{M})^{\frac1r}$. By Theorem \ref{thm:joint-discrete-coefficient-error}, we may take $e_\alpha\sim \delta^\beta$. Then $\rho\sim\delta^\beta$, and the parameter choices $\gamma\sim\delta^{2\beta}$, $h\sim\delta^{\frac\beta2}$  and $\tau=O(\delta^\beta)$   
yield $\eta_*=O(\delta^\beta)$. It follows from the estimate \eqref{eq:joint-discrete-coefficient-error} that
\begin{align*}
	\|q_h^*-q^\dagger\|_{L^2(\Omega)} \le c_T\delta^{\frac\beta4} =c_T\delta^{\frac14-\frac{5}{8r}}.
\end{align*}
In particular.  for a piecewise quadratic reconstruction ($p=2$) with $s\ge3$, we have $r=3$ and $\beta=\frac16$, and thus $\|q_h^*-q^\dagger\|_{L^2(\Omega)} \le c_T \delta^{\frac{1}{24}}$.
\end{remark}

\section{Conclusion}
In this work, we have investigated the stability and reconstruction for the inverse problem of identifying the diffusion coefficient in a source-free time-fractional diffusion model from a single terminal observation. Unlike classical diffusion, fractional evolution does not exhibit first-mode dominance at large times, since all spectral modes decay at the same algebraic rate, and moreover, the absence of a source means that the terminal state approaches zero. We address these challenges through a careful splitting that involves an elliptic leading term and a coefficient perturbation estimate for the evolution remainder.

The first main contribution is a conditional stability estimate under asymmetric regularity assumptions: stronger regularity is imposed only on the true coefficient $q^\dag$, and the comparison coefficient $q$ requires only $H^1$ regularity. This distinction provides a suitable framework for numerical analysis, since recovered coefficient need not share the smoothness of the true coefficient. The key is a perturbation estimate for the remainder that does not involve higher-order norms of the comparison coefficient. Furthermore, we prove that the time threshold is independent of its gradient, and thus the observation time can be fixed independently of the mesh size and regularization parameter. These properties provide the basis for the subsequent regularization and discretization error analysis.

The second main contribution is a boundary identity that separates the fractional order from the unknown interior coefficient. When the coefficient is known in a neighborhood of the boundary patch and the initial datum has a nonvanishing trace on the boundary portion, the additional regularity of the remainder makes its contribution to the relevant boundary trace vanish. This identity enables identifying the fractional order before reconstructing the diffusion coefficient, laying the basis for their sequential recovery from the same terminal observation. Finally, the regularized reconstruction and error analysis extend the continuous stability results to a quantitative analysis of the numerical recovery. 

Together, these results provide a theoretical framework for parameter identification from a single decaying terminal state, from conditional stability and boundary identification to regularization and discretization error estimates.

\bibliographystyle{abbrv}
\bibliography{reference}

@article {KianLiLiu:2021,
    AUTHOR = {Kian, Yavar and Li, Zhiyuan and Liu, Yikan and Yamamoto,
              Masahiro},
     TITLE = {The uniqueness of inverse problems for a fractional equation
              with a single measurement},
   JOURNAL = {Math. Ann.},
  FJOURNAL = {Mathematische Annalen},
    VOLUME = {380},
      YEAR = {2021},
    NUMBER = {3-4},
     PAGES = {1465--1495},
      ISSN = {0025-5831,1432-1807},
   MRCLASS = {35R30 (35R11 58J99)},
  MRNUMBER = {4297191},
       DOI = {10.1007/s00208-020-02027-z},
       URL = {https://doi-org.easyaccess1.lib.cuhk.edu.hk/10.1007/s00208-020-02027-z},
}

@article {Boulakia:2026,
    AUTHOR = {Boulakia, Muriel and Burman, Erik and Fern\'andez, Miguel A.},
     TITLE = {Primal-dual finite element methods for the direct
              reconstruction of material coefficients in stationary elliptic
              problems},
   JOURNAL = {Inverse Problems},
  FJOURNAL = {Inverse Problems. An International Journal on the Theory and
              Practice of Inverse Problems, Inverse Methods and Computerized
              Inversion of Data},
    VOLUME = {42},
      YEAR = {2026},
    NUMBER = {6},
     PAGES = {065004, 23 pp.},
      ISSN = {0266-5611,1361-6420},
   MRCLASS = {65N21 (35J15 35R30 74G75 74S05)},
  MRNUMBER = {5083983},
       DOI = {10.1088/1361-6420/ae72be},
       URL = {https://doi-org.easyaccess1.lib.cuhk.edu.hk/10.1088/1361-6420/ae72be},
}

@article {Alessandiri:1986,
    AUTHOR = {Alessandrini, Giovanni},
     TITLE = {An identification problem for an elliptic equation in two
              variables},
   JOURNAL = {Ann. Mat. Pura Appl. (4)},
  FJOURNAL = {Annali di Matematica Pura ed Applicata. Serie Quarta},
    VOLUME = {145},
      YEAR = {1986},
     PAGES = {265--295},
      ISSN = {0003-4622},
   MRCLASS = {35R30 (35J25)},
  MRNUMBER = {886713},
MRREVIEWER = {S.\ D\"ummel},
       DOI = {10.1007/BF01790543},
       URL = {https://doi-org.easyaccess1.lib.cuhk.edu.hk/10.1007/BF01790543},
}

@article {Falk:1983,
    AUTHOR = {Falk, Richard S.},
     TITLE = {Error estimates for the numerical identification of a variable
              coefficient},
   JOURNAL = {Math. Comp.},
  FJOURNAL = {Mathematics of Computation},
    VOLUME = {40},
      YEAR = {1983},
    NUMBER = {162},
     PAGES = {537--546},
      ISSN = {0025-5718,1088-6842},
   MRCLASS = {65N99 (65N30)},
  MRNUMBER = {689469},
MRREVIEWER = {E.\ G.\ D\cprime yakonov},
       DOI = {10.2307/2007529},
       URL = {https://doi-org.easyaccess1.lib.cuhk.edu.hk/10.2307/2007529},
}

@article {Richter:1981,
    AUTHOR = {Richter, Gerard R.},
     TITLE = {An inverse problem for the steady state diffusion equation},
   JOURNAL = {SIAM J. Appl. Math.},
  FJOURNAL = {SIAM Journal on Applied Mathematics},
    VOLUME = {41},
      YEAR = {1981},
    NUMBER = {2},
     PAGES = {210--221},
      ISSN = {0036-1399},
   MRCLASS = {35R30 (76R99)},
  MRNUMBER = {628945},
MRREVIEWER = {S.\ D\"ummel},
       DOI = {10.1137/0141016},
       URL = {https://doi-org.easyaccess1.lib.cuhk.edu.hk/10.1137/0141016},
}

@article{cheng2009uniqueness,
	title={Uniqueness in an inverse problem for a one-dimensional fractional diffusion equation},
	author={Cheng, Jin and Nakagawa, Junichi and Yamamoto, Masahiro and Yamazaki, Tomohiro},
	journal={Inverse problems},
	volume={25},
	number={11},
	pages={115002, 16 pp.},
	year={2009},
        doi = {10.1088/0266-5611/25/11/115002},
        url = {https://doi.org/10.1088/0266-5611/25/11/115002},
	publisher={IOP Publishing}
}

@article{LinLin:2026,
  title={Global recovery of a scalar diffusivity from a single final state},
  author={Lin, Ching-Lung and Lin, Yi-Hsuan},
  journal={preprint, arXiv:2610.32063},
  year={2026}
}

@book{ErnGuermond2021,
  author    = {Ern, Alexandre and Guermond, Jean-Luc},
  title     = {{Finite Elements I: Approximation and Interpolation}},
  publisher = {Springer, Cham},
  year      = {2021},
  doi       = {10.1007/978-3-030-56341-7}
}

@article {LiYangZhou:2024,
    AUTHOR = {Li, Buyang and Yang, Zongze and Zhou, Zhi},
     TITLE = {High-order splitting finite element methods for the
              subdiffusion equation with limited smoothing property},
   JOURNAL = {Math. Comp.},
  FJOURNAL = {Mathematics of Computation},
    VOLUME = {93},
      YEAR = {2024},
    NUMBER = {350},
     PAGES = {2557--2586},
      ISSN = {0025-5718,1088-6842},
   MRCLASS = {65M30 (65M12 65M15)},
  MRNUMBER = {4780338},
MRREVIEWER = {Mohadese\ Ramezani},
       DOI = {10.1090/mcom/3944},
       URL = {https://doi.org/10.1090/mcom/3944},
}

@article {JinLuQuanZhou:2025,
    AUTHOR = {Jin, Bangti and Lu, Xiliang and Quan, Qimeng and Zhou, Zhi},
     TITLE = {Numerical recovery of the diffusion coefficient in diffusion
              equations from terminal measurement},
   JOURNAL = {Numer. Math.},
  FJOURNAL = {Numerische Mathematik},
    VOLUME = {157},
      YEAR = {2025},
    NUMBER = {6},
     PAGES = {2323--2355},
      ISSN = {0029-599X,0945-3245},
   MRCLASS = {65M32 (65M15 65M30 65M60)},
  MRNUMBER = {4984394},
       DOI = {10.1007/s00211-025-01495-2},
       URL = {https://doi.org/10.1007/s00211-025-01495-2},
}

@book{JinZhou:2023book,
   author = {Jin, Bangti and Zhou, Zhi},
    title = {{Numerical Treatment and Analysis of Time-Fractional Evolution Equations}},
     year = {2023},
  publisher = {Springer, Cham},
}

@article{nigmatullin1986realization,
	title={The realization of the generalized transfer equation in a medium with fractal geometry},
	author={Nigmatullin, R R},
	journal={Physica Status Solidi (b)},
	volume={133},
	number={1},
	pages={425--430},
	year={1986},
         url = {https://doi.org/10.1002/pssb.2221330150},
         doi = {10.1002/pssb.2221330150}, 
	publisher={Wiley Online Library}
}

@article{lubich1986discretized,
    AUTHOR = {Lubich, Ch.},
     TITLE = {Discretized fractional calculus},
   JOURNAL = {SIAM J. Math. Anal.},
  FJOURNAL = {SIAM Journal on Mathematical Analysis},
    VOLUME = {17},
      YEAR = {1986},
    NUMBER = {3},
     PAGES = {704--719},
      ISSN = {0036-1410},
   MRCLASS = {26A33 (41A55 65D30)},
  MRNUMBER = {838249},
MRREVIEWER = {R. Bittner},
       DOI = {10.1137/0517050},
       URL = {https://doi.org/10.1137/0517050},
}

@article{giona1992fractional,
	title={Fractional diffusion equation for transport phenomena in random media},
	author={Giona, Massimiliano and Roman, H Eduardo},
	journal={Phys. A: Stat. Mech. Appl.},
	volume={185},
	number={1-4},
	pages={87--97},
	year={1992},
         url = {https://doi.org/10.1016/0378-4371(92)90441-R},
         doi = {https://doi.org/10.1016/0378-4371(92)90441-R},
	publisher={Elsevier}
}

@book{engl1996regularization,
    AUTHOR = {Engl, Heinz W. and Hanke, Martin and Neubauer, Andreas},
     TITLE = {{Regularization of Inverse Problems}},
 PUBLISHER = {Kluwer Academic Publishers Group, Dordrecht},
      YEAR = {1996},
     PAGES = {viii+321},
      ISBN = {0-7923-4157-0},
   MRCLASS = {65J20 (35R30 47A50 47H17)},
  MRNUMBER = {1408680},
MRREVIEWER = {Ulrich Tautenhahn},
}

@article{metzler2000random,
	title={The random walk's guide to anomalous diffusion: a fractional dynamics approach},
	author={Metzler, Ralf and Klafter, Joseph},
	journal={Phys. Rep.},
	volume={339},
	number={1},
	pages={1--77},
	year={2000},
         url = {https://doi.org/10.1016/S0370-1573(00)00070-3},
         doi = {10.1016/S0370-1573(00)00070-3},
	publisher={Elsevier}
}

@article{adams1992field,
	title={Field study of dispersion in a heterogeneous aquifer: 2. Spatial moments analysis},
	author={Adams, E Eric and Gelhar, Lynn W},
	journal={Water Resources Res.},
	volume={28},
	number={12},
	pages={3293--3307},
	year={1992},
          doi = {10.1029/92WR01757}, 
    url ={https://doi.org/10.1029/92WR01757},
	publisher={Wiley Online Library}
}

@article{hatano1998dispersive,
	title={Dispersive transport of ions in column experiments: An explanation of long-tailed profiles},
	author={Hatano, Yuko and Hatano, Naomichi},
	journal={Water Resources Res.},
	volume={34},
	number={5},
	pages={1027--1033},
	year={1998},
         url = {10.1029/98WR00214},
         doi = { https://doi.org/10.1029/98WR00214},
	publisher={Wiley Online Library}
}

@book{kilbas2006theory,
    AUTHOR = {Kilbas, Anatoly A. and Srivastava, Hari M. and Trujillo, Juan
              J.},
     TITLE = {{Theory and Applications of Fractional Differential Equations}},
 PUBLISHER = {Elsevier Science B.V., Amsterdam},
      YEAR = {2006},
     PAGES = {xvi+523},
      ISBN = {978-0-444-51832-3; 0-444-51832-0},
   MRCLASS = {34-02 (26A33 33C90 34A99 35-02 45-02)},
  MRNUMBER = {2218073},
MRREVIEWER = {B. S. Rubin},
}

@article {Triki:2021,
    AUTHOR = {Triki, Faouzi},
     TITLE = {Coefficient identification in parabolic equations with final
              data},
   JOURNAL = {J. Math. Pures Appl. (9)},
  FJOURNAL = {Journal de Math\'{e}matiques Pures et Appliqu\'{e}es. Neuvi\`eme S\'{e}rie},
    VOLUME = {148},
      YEAR = {2021},
     PAGES = {342--359},
      ISSN = {0021-7824},
   MRCLASS = {35R30 (35K15 35K20)},
  MRNUMBER = {4223356},
MRREVIEWER = {Enno Pais},
       DOI = {10.1016/j.matpur.2021.02.004},
       URL = {https://doi.org/10.1016/j.matpur.2021.02.004},
}

@article{Bonito:2017,
    AUTHOR = {Bonito, Andrea and Cohen, Albert and DeVore, Ronald and
              Petrova, Guergana and Welper, Gerrit},
     TITLE = {Diffusion coefficients estimation for elliptic partial
              differential equations},
   JOURNAL = {SIAM J. Math. Anal.},
  FJOURNAL = {SIAM Journal on Mathematical Analysis},
    VOLUME = {49},
      YEAR = {2017},
    NUMBER = {2},
     PAGES = {1570--1592},
      ISSN = {0036-1410},
   MRCLASS = {35R30 (35J57)},
  MRNUMBER = {3639575},
MRREVIEWER = {Aleksander Denisiuk},
       DOI = {10.1137/16M1094476},
       URL = {https://doi.org/10.1137/16M1094476},
}

@book {GilbargTrudinger:2001,
    AUTHOR = {Gilbarg, David and Trudinger, Neil S.},
     TITLE = {{Elliptic Partial Differential Equations of Second Order}},
 PUBLISHER = {Springer-Verlag, Berlin},
      YEAR = {2001},
     PAGES = {xiv+517},
      ISBN = {3-540-41160-7},
   MRCLASS = {35-02 (35Jxx)},
  MRNUMBER = {1814364},
}

@article{JinZhou:2021SINUM,
    AUTHOR = {Jin, Bangti and Zhou, Zhi},
     TITLE = {Error analysis of finite element approximations of diffusion
              coefficient identification for elliptic and parabolic
              problems},
   JOURNAL = {SIAM J. Numer. Anal.},
  FJOURNAL = {SIAM Journal on Numerical Analysis},
    VOLUME = {59},
      YEAR = {2021},
    NUMBER = {1},
     PAGES = {119--142},
      ISSN = {0036-1429},
   MRCLASS = {65N21 (35R11 65M15 65M32 65N15)},
  MRNUMBER = {4196420},
       DOI = {10.1137/20M134383X},
       URL = {https://doi.org/10.1137/20M134383X},
}

@book {ItoJin:2015,
	AUTHOR = {Ito, Kazufumi and Jin, Bangti},
	TITLE = {Inverse {P}roblems: Tikhonov {T}heory and {A}lgorithms},
	PUBLISHER = {World Scientific Publishing Co. Pte. Ltd., Hackensack, NJ},
	YEAR = {2015},
	PAGES = {x+318},
	ISBN = {978-981-4596-19-0},
	MRCLASS = {65-02 (49N45 65J22 65M32 65N21)},
	MRNUMBER = {3244283},
	MRREVIEWER = {Sudhir H. Kulkarni},
}

@article {Alessandrini:2020,
    AUTHOR = {Alessandrini, Giovanni},
     TITLE = {A small collection of open problems},
   JOURNAL = {Rend. Istit. Mat. Univ. Trieste},
  FJOURNAL = {Rendiconti dell'Istituto di Matematica dell'Universit\`a di
              Trieste. An International Journal of Mathematics},
    VOLUME = {52},
      YEAR = {2020},
     PAGES = {591--600},
      ISSN = {0049-4704},
   MRCLASS = {35K10},
  MRNUMBER = {4207653},
       DOI = {10.1007/s11118-018-9751-1},
       URL = {https://doi.org/10.1007/s11118-018-9751-1},
}

@article {AlessandriniVessella:1985,
    AUTHOR = {Alessandrini, Giovanni and Vessella, Sergio},
     TITLE = {Error estimates in an identification problem for a parabolic
              equation},
   JOURNAL = {Boll. Un. Mat. Ital. C (6)},
  FJOURNAL = {Unione Matematica Italiana. Bollettino. C. Serie VI. Analisi
              Funzionale e Applicazioni},
    VOLUME = {4},
      YEAR = {1985},
    NUMBER = {1},
     PAGES = {183--203},
   MRCLASS = {35R30 (35K05)},
  MRNUMBER = {805213},
MRREVIEWER = {Marian Majchrowski},
}

@article {Isakov:1991,
    AUTHOR = {Isakov, Victor},
     TITLE = {Inverse parabolic problems with the final overdetermination},
   JOURNAL = {Comm. Pure Appl. Math.},
  FJOURNAL = {Communications on Pure and Applied Mathematics},
    VOLUME = {44},
      YEAR = {1991},
    NUMBER = {2},
     PAGES = {185--209},
      ISSN = {0010-3640},
   MRCLASS = {35R30 (35K10 35N10)},
  MRNUMBER = {1085828},
MRREVIEWER = {Dang Dinh Ang},
       DOI = {10.1002/cpa.3160440203},
       URL = {https://doi.org/10.1002/cpa.3160440203},
}

@book {Kubica:book2020,
    AUTHOR = {Kubica, Adam and Ryszewska, Katarzyna and Yamamoto, Masahiro},
     TITLE = {{Time-Fractional Differential Equations---a Theoretical
    Introduction}},
 PUBLISHER = {Springer, Singapore},
      YEAR = {2020},
     PAGES = {x+134},
      ISBN = {978-981-15-9066-5; 978-981-15-9065-8},
   MRCLASS = {35R11 (34-01 34A08 35-01)},
  MRNUMBER = {4200127},
MRREVIEWER = {Paulo\ M.\ Carvalho-Neto},
       DOI = {10.1007/978-981-15-9066-5},
       URL = {https://doi.org/10.1007/978-981-15-9066-5},
}

@book {KaltenbacherRundell:book2023,
    AUTHOR = {Kaltenbacher, Barbara and Rundell, William},
     TITLE = {{Inverse Problems for Fractional Partial Differential
    Equations}},
 PUBLISHER = {AMS, Providence, RI},
      YEAR = {2023},
     PAGES = {xiii+505},
      ISBN = {[9781470472450]; [9781470472764]},
   MRCLASS = {35-02 (26A33 35K57 35R11 35R30 47F05 60K50 65M32)},
  MRNUMBER = {4604107},
       DOI = {10.1090/gsm/230},
       URL = {https://doi.org/10.1090/gsm/230},
}

@book {Jin:book2021,
    AUTHOR = {Jin, Bangti},
     TITLE = {{Fractional Differential Equations --- An Approach via Fractional Derivatives}},
 PUBLISHER = {Springer, Cham},
      YEAR = {2021},
     PAGES = {xiv+368},
      ISBN = {978-3-030-76042-7; 978-3-030-76043-4},
   MRCLASS = {35R11 (26A33 34A08)},
  MRNUMBER = {4290515},
       DOI = {10.1007/978-3-030-76043-4},
       URL = {https://doi.org/10.1007/978-3-030-76043-4},
}

@article {JinLiZhou:2017,
    AUTHOR = {Jin, Bangti and Li, Buyang and Zhou, Zhi},
     TITLE = {Correction of high-order {BDF} convolution quadrature for
              fractional evolution equations},
   JOURNAL = {SIAM J. Sci. Comput.},
  FJOURNAL = {SIAM Journal on Scientific Computing},
    VOLUME = {39},
      YEAR = {2017},
    NUMBER = {6},
     PAGES = {A3129--A3152},
      ISSN = {1064-8275,1095-7197},
   MRCLASS = {65M60 (35R11 65M12 65M15)},
  MRNUMBER = {3738850},
MRREVIEWER = {V.\ L.\ Makarov},
       DOI = {10.1137/17M1118816},
       URL = {https://doi.org/10.1137/17M1118816},
}

@article {Stynes:2017,
    AUTHOR = {Stynes, Martin and O'Riordan, Eugene and Gracia, Jos\'{e}
              Luis},
     TITLE = {Error analysis of a finite difference method on graded meshes
              for a time-fractional diffusion equation},
   JOURNAL = {SIAM J. Numer. Anal.},
  FJOURNAL = {SIAM Journal on Numerical Analysis},
    VOLUME = {55},
      YEAR = {2017},
    NUMBER = {2},
     PAGES = {1057--1079},
      ISSN = {0036-1429,1095-7170},
   MRCLASS = {65M06 (35B65 35R11 65M15)},
  MRNUMBER = {3639581},
MRREVIEWER = {Ali\ R.\ Soheili},
       DOI = {10.1137/16M1082329},
       URL = {https://doi.org/10.1137/16M1082329},
}

@article {JinLazarovZhou:2013,
    AUTHOR = {Jin, Bangti and Lazarov, Raytcho and Zhou, Zhi},
     TITLE = {Error estimates for a semidiscrete finite element method for
              fractional order parabolic equations},
   JOURNAL = {SIAM J. Numer. Anal.},
  FJOURNAL = {SIAM Journal on Numerical Analysis},
    VOLUME = {51},
      YEAR = {2013},
    NUMBER = {1},
     PAGES = {445--466},
      ISSN = {0036-1429,1095-7170},
   MRCLASS = {65M60 (26A33 35R11 45K05 65M15)},
  MRNUMBER = {3033018},
MRREVIEWER = {Akhtar\ A.\ Khan},
       DOI = {10.1137/120873984},
       URL = {https://doi.org/10.1137/120873984},
}

@article {LiSun:2017,
    AUTHOR = {Li, Buyang and Sun, Weiwei},
     TITLE = {Maximal {$L^p$} analysis of finite element solutions for
              parabolic equations with nonsmooth coefficients in convex
              polyhedra},
   JOURNAL = {Math. Comp.},
  FJOURNAL = {Mathematics of Computation},
    VOLUME = {86},
      YEAR = {2017},
    NUMBER = {305},
     PAGES = {1071--1102},
      ISSN = {0025-5718,1088-6842},
   MRCLASS = {65M60 (35K20 65M12)},
  MRNUMBER = {3614012},
MRREVIEWER = {B\"{u}lent\ Karas\"{o}zen},
       DOI = {10.1090/mcom/3133},
       URL = {https://doi.org/10.1090/mcom/3133},
}

@article {VergaraZacher:2015,
    AUTHOR = {Vergara, Vicente and Zacher, Rico},
     TITLE = {Optimal decay estimates for time-fractional and other nonlocal
              subdiffusion equations via energy methods},
   JOURNAL = {SIAM J. Math. Anal.},
  FJOURNAL = {SIAM Journal on Mathematical Analysis},
    VOLUME = {47},
      YEAR = {2015},
    NUMBER = {1},
     PAGES = {210--239},
      ISSN = {0036-1410,1095-7154},
   MRCLASS = {45K05 (35B40 35K59 35K92 35R11 47G20)},
  MRNUMBER = {3296607},
MRREVIEWER = {Dmitry\ E.\ Pelinovsky},
       DOI = {10.1137/130941900},
       URL = {https://doi.org/10.1137/130941900},
}

@article {Allen:2016,
    AUTHOR = {Allen, Mark and Caffarelli, Luis and Vasseur, Alexis},
     TITLE = {A parabolic problem with a fractional time derivative},
   JOURNAL = {Arch. Ration. Mech. Anal.},
  FJOURNAL = {Archive for Rational Mechanics and Analysis},
    VOLUME = {221},
      YEAR = {2016},
    NUMBER = {2},
     PAGES = {603--630},
      ISSN = {0003-9527,1432-0673},
   MRCLASS = {35R11 (26A33 35A01 35A02 35B65)},
  MRNUMBER = {3488533},
MRREVIEWER = {C\'{e}sar\ Enrique\ Torres Ledesma},
       DOI = {10.1007/s00205-016-0969-z},
       URL = {https://doi.org/10.1007/s00205-016-0969-z},
}

@article {Sakamoto:2011,
    AUTHOR = {Sakamoto, Kenichi and Yamamoto, Masahiro},
     TITLE = {Initial value/boundary value problems for fractional
              diffusion-wave equations and applications to some inverse
              problems},
   JOURNAL = {J. Math. Anal. Appl.},
  FJOURNAL = {Journal of Mathematical Analysis and Applications},
    VOLUME = {382},
      YEAR = {2011},
    NUMBER = {1},
     PAGES = {426--447},
      ISSN = {0022-247X,1096-0813},
   MRCLASS = {35R11 (35B30 35B35 35B40 35B45 35R30)},
  MRNUMBER = {2805524},
MRREVIEWER = {Giulio\ Ciraolo},
       DOI = {10.1016/j.jmaa.2011.04.058},
       URL = {https://doi.org/10.1016/j.jmaa.2011.04.058},
}

@article{ouhabaz1995gaussian,
    AUTHOR = {Ouhabaz, El-Maati},
     TITLE = {Gaussian estimates and holomorphy of semigroups},
   JOURNAL = {Proc. Amer. Math. Soc.},
  FJOURNAL = {Proceedings of the American Mathematical Society},
    VOLUME = {123},
      YEAR = {1995},
    NUMBER = {5},
     PAGES = {1465--1474},
      ISSN = {0002-9939},
   MRCLASS = {47D06 (47F05 47N20)},
  MRNUMBER = {1232142},
MRREVIEWER = {J. A. van Casteren},
       DOI = {10.2307/2161136},
       URL = {https://doi.org/10.2307/2161136},
}
\end{document}